\documentclass{amsart}
\usepackage{latexsym,epsfig,fullpage,amsfonts,amssymb,amsmath,amscd,epic,graphics}
\newtheorem{theorem}{Theorem}[section]
\newtheorem{lemma}[theorem]{Lemma}
\newtheorem{conjecture}[theorem]{Conjecture}
\newtheorem{proposition}[theorem]{Proposition}
\theoremstyle{definition}
\newtheorem{definition}[theorem]{Definition}

\theoremstyle{remark}
\newtheorem{remark}[theorem]{Remark}

\numberwithin{equation}{section}

\newcommand{\Sig}{{\Sigma}}
\newcommand{\PP}{{\mathbb P}}
\newcommand{\op}{{\overline p}}
\newcommand{\C}{{\mathbb C}}
\newcommand{\K}{{\mathbb K}}
\newcommand{\val}{{\operatorname{Val}}}
\newcommand{\red}{{\operatorname{red}}}
\newcommand{\Int}{{\operatorname{Int}}}
\newcommand{\const}{{\operatorname{const}}}
\newcommand{\eps}{{\varepsilon}}
\newcommand{\R}{{\mathbb R}}
\newcommand{\bw}{{\boldsymbol w}}
\newcommand{\Q}{{\mathbb Q}}
\newcommand{\Del}{{\Delta}}\newcommand{\del}{{\delta}}
\newcommand{\lam}{{\lambda}}\newcommand{\Lam}{{\Lambda}}
\newcommand{\Gam}{{\Gamma}}
\newcommand{\Tor}{{\operatorname{Tor}}}
\newcommand{\Aff}{{\operatorname{Aff}}}
\newcommand{\Z}{{\mathbb Z}}
\newcommand{\conv}{{\operatorname{Conv}}}
\newcommand{\gam}{{\gamma}}
\newcommand{\sig}{{\sigma}}
\newcommand{\Iso}{{\operatorname{Iso}}}
\newcommand{\rk}{{\operatorname{rk}}}
\newcommand{\bp}{{\boldsymbol p}}\newcommand{\by}{{\boldsymbol y}}
\newcommand{\bx}{{\boldsymbol x}}
\newcommand{\conj}{{\operatorname{Conj}}}
\newcommand{\alp}{{\alpha}}
\newcommand{\bet}{{\beta}}
\newcommand{\Sing}{{\operatorname{Sing}}}

\begin{document}

\title{A tropical calculation of the Welschinger invariants of real toric Del Pezzo surfaces}
\author{Eugenii Shustin}
\address{School of Mathematical Sciences, Tel Aviv University, Ramat Aviv, 69978 Tel Aviv, Israel}
\email{shustin@post.tau.ac.il}
\thanks{Part of this work was done during the
author's stay at Universit\"at Kaiserslautern, supported by the
Hermann-Minkowski Minerva Center for Geometry at Tel Aviv
University, and during the author's stay at the Mathematical
Science Research Institute, Berkeley. The author is very grateful
to the Hermann-Minkowski Minerva Center for its support, and to
Universit\"at Kaiserslautern and MSRI for the hospitality and
excellent work conditions.} 

\date{}


\begin{abstract}
The Welschinger invariants of real rational algebraic surfaces are
natural analogues of the genus zero Gromov-Witten invariants. We
establish a tropical formula to calculate the Welschinger
invariants of real toric Del Pezzo surfaces for any
conjugation-invariant configuration of points. The formula
expresses the Welschinger invariants via the total multiplicity of
certain tropical curves (non-Archimedean amoebas) passing through
generic configurations of points, and then via the total
multiplicity of some lattice paths in the convex lattice polygon
associated with a given surface. We also present the results of
computation of Welschinger invariants, obtained jointly with I.
Itenberg and V. Kharlamov.
\end{abstract}

\maketitle

\section*{Introduction}

One of the important questions of real enumerative geometry is:
for a real algebraic surface $\Sig$, how many real rational curves
in an ample linear system $|D|$ pass through a
conjugation-invariant configuration of $-K_\Sig D-1$ distinct
generic points in $\Sig$? Similarly to the complex case, where the
answer is given by Gromov-Witten invariants, the recently
discovered Welschinger invariants \cite{Wel,Wel1} appear to be an
ultimate tool to handle the question over the reals. In
particular, they led to non-trivial positive lower bounds for the
numbers in question, provided that the configuration consisted of
only real points \cite{IKS}.

So far no closed or recursive formula is found for the Welschinger
invariants, and tropical enumerative geometry \cite{M2,M3,ShP}
provides the only known approach to compute them. In the present
paper we express the Welschinger invariants for configurations of
real and imaginary conjugate points on real toric Del Pezzo
surfaces via the number of certain subdivisions of the
corresponding convex lattice polygons.

\medskip
\noindent {\bf Welschinger invariants.} Let $\Sig$ be $\PP^2$, or
the hyperboloid $(\PP^1)^2$, or the plane $\PP^2_k$ blown up at
$k$ generic real points, equipped with the standard real
structure. Let $D$ be a real ample divisor on $\Sig$, and let the
non-negative integers $r',r''$ satisfy
\begin{equation}r'+2r''=-K_\Sig D-1\
.\label{enn36}\end{equation} Denote by $\Omega_{r''}(\Sig,D)$ the
set of configurations of $-K_\Sigma D - 1$ distinct points of
$\Sig$ such that $r'$ of them are real and the rest form $r''$
pairs of imaginary conjugate points. The Welschinger number
$W_{r''}(\Sig,D)$ is the sum of weights of all the real rational
curves in $|D|$, passing through a generic
configuration\footnote{The generality means here that all the
complex rational curves through the configuration are nodal
irreducible, and their number is equal to the corresponding
Gromov-Witten invariant.} $\op\in\Omega_{r''}(\Sig,D)$, where the
weight of a real rational curve $C$ is $1$ if it has an even
number of real solitary nodes, and is $-1$ otherwise. The surfaces
$\Sig$ as above are the only toric surfaces, whose complex
structure determines a symplectic structure which is generic in
its deformation class, and thus, by Welschinger's theorem
\cite{Wel,Wel1}, $W_{r''}(\Sig,D)$ does not depend on the choice
of a generic element $\op\in\Omega_{r''}(\Sig,D)$ (a simple proof
of the latter independence in the algebraic setting is found in
\cite{IKS1}). The importance of the Welschinger invariant comes
from an immediate inequality
\begin{equation}|W_{r''}(\Sig,D)|\le R_{\Sig,D}(\op)\le N_{\Sig,D}\
,\label{enn200}\end{equation} where $R_{\Sig,D}(\op)$ is the
number of real rational curves in $|D|$ passing through a generic
configuration $\op\in\Omega_{r''}(\Sig,D)$, and $N_{\Sig,D}$ is
the number of complex rational curves in $|D|$, passing through
generic $-K_\Sigma D - 1$ points in $\Sig$. We should notice that
even the existence of a positive lower bound for $R_{\Sig,D}(\op)$
was not reached by any other method, beginning with the case of
plane quartics.

\medskip
\noindent {\bf A tropical calculation of the Welschinger
invariant.} Our approach to calculating the Welschinger invariant
is quite similar to that in \cite{IKS}, and it heavily relies on
the enumerative tropical algebraic geometry developed in
\cite{M2,M3,ShP}. More precisely, we replace the complex field
$\C$ by the field $\K=\bigcup_{m\ge 1}\C\{\{t^{1/m}\}\}$ of the
complex locally convergent Puiseux series equipped with the
standard complex conjugation and with a non-Archimedean valuation
$$\val:\K^*\to\R,\quad\val\left(\sum_ka_kt^k\right)=-\min\{k\ :\ a_k\ne 0\}\
.$$ A rational curve over $\K_\R$ passing through a generic
configuration $\op\in\Omega_{r''}(\Sig_\K,D)$ is represented as an
equisingular family of real rational curves in $\Sig$ over the
punctured disc, and its limit at the disc center determines a
tropical curve in $\R^2$ (called a {\it real rational tropical
curve}), which passes through the configuration
$\val(\op)\subset\R^2$, coordinate-wise $\val$-projection of
$\op$.

Our first main result is Theorem \ref{tnn1} (section
\ref{secn14}), which precisely describes the real rational
tropical curves, passing through generic configurations of points
in $\Q^2$, and, for any such real rational tropical curve $A$,
determines the contribution to the Welschinger invariant of all
the real algebraic curves, projecting to $A$. The proof is based
on the techniques of \cite{ShP}.

In the case of configurations, consisting of only real points,
this result has been obtained in \cite{M3,ShP}. We notice that an
extension of the tropical formula to the case of configurations of
real and imaginary conjugate points, requires some further
development of the techniques of tropical enumerative geometry,
which we present in this paper. The principal difficulty in the
latter case is that the configuration $\val(\op)$ contains fewer
points than $\op$, since a pair of conjugate points projects to
the same point in $\R^2$.

The second main result is Theorem \ref{cnn1} (section
\ref{secn10}), which reduces the count of real rational tropical
curves, passing through a generic configuration in $\Q^2$, to the
count of the total weight of certain lattice paths in the convex
lattice polygon corresponding to the divisor $D$. The weight of a
lattice path is the sum of Welschinger numbers of certain
subdivisions of the given polygon into convex lattice subpolygons,
which can be produced from the lattice path in a finite
combinatorial algorithm (section \ref{secn10}). Here we follow the
Mikhalkin's idea to place the configuration $\val(\op)\subset\R^2$
on a straight line.

\medskip
\noindent{\bf Applications.} The tropical formula turns the
computation of Welschinger invariants into a purely combinatorial
problem on geometry of lattice paths and lattice subdivisions of
convex lattice polygons. We present here some results, obtained in
this way jointly with I. Itenberg and V. Kharlamov (see
\cite{IKS2} for a detailed presentation and proofs). These results
concern the positivity of Welschinger invariants, their monotone
behavior with respect to the number of imaginary points in
configurations, and the asymptotics with respect to the growing
degree of the divisor $D$. We formulate a few natural conjectures.

\medskip
{\it A. The positivity and monotonicity of the Welschinger
invariant.} The positivity of $W_0(\Sig,D)$ for all real toric Del
Pezzo surfaces and ample divisors $D$ was shown in \cite{IKS}. For
$r''\ne 0$, the invariants $W_{r''}(\Sig,D)$ can vanish (see, for
example, computation of $W_{r''}(\PP^2,3L)$, $L$ being a line in
the plane, in section \ref{secn12}).

\begin{theorem}\label{tn5} {\rm (\cite{IKS2})} (1) The invariants
$W_{r''}((\PP^1)^2,D)$ are positive for all ample divisors $D$ and
all $r''\ge 0$, satisfying (\ref{enn36}).

(2) The following inequalities hold: \begin{enumerate}\item[(i)]
$W_0(\Sig,D)>W_1(\Sig,D)>W_2(\Sig,D)>W_3(\Sig,D)$, for
$\Sig=\PP^2$, $(\PP^1)^2$, $\PP^2_1$, or $\PP^2_2$, with $D$ an
ample divisor such that $D(D+K_\Sig)\ge 0$; \item[(ii)]
$W_3(\PP^2,D)>0$, and $W_2(\Sig,D)>0$, for $\Sig=\PP^2_k$,
$k=1,2,3$, and an ample divisor $D$.\end{enumerate}
\end{theorem}

\begin{conjecture} {\rm (\cite{IKS2})}
For a real unnodal Del Pezzo surface $\Sig$ and any ample divisor
$D$ on $\Sig$, the Welschinger invariants $W_{r''}(\Sig,D)$ are
positive as $0\le r''<[(-K_\Sig D-1)/2]$, and are non-negative for
$r''=[(-K_\Sig D-1)/2]$. Furthermore, they satisfy the
monotonicity relation
$$W_{r''}(\Sig,D)\ge W_{r''+1}(\Sig,D),\quad 0\le r''<\left[\frac{-K_\Sig
D-1}{2}\right] \ .$$
\end{conjecture}

We notice that the monotonicity and non-negativity of Welschinger
invariants are closely related, since, by \cite{Wel}, Theorem 2.2,
the first difference of the function $r''\mapsto W_{r''}(\Sig,D)$
is twice the Welschinger invariant for the surface $\Sig$ blown up
at one real point.

\medskip
{\it B. The asymptotics of the Welschinger invariants.}

\begin{theorem}\label{tn4} {\rm (\cite{IKS2})}
The Welschinger invariants of the plane satisfy the relation
\begin{equation}\lim_{n\to\infty}\frac{\log
W_{r''}(\PP^2,nL)}{n\log n}=\lim_{n\to\infty}\frac{\log
N_{\PP^2,nL}}{n\log n}=3,\quad 0\le r''\le 3\
.\label{enn60}\end{equation} For $\Sig=(\PP^1)^2$, or $\PP^2_k$,
$1\le k\le 3$, it holds that
\begin{equation}\lim_{d\to\infty}\frac{\log
W_{r''}(\Sig,nD)}{n\log n}=\lim_{n\to\infty}\frac{\log
N_{\Sig,nD}}{n\log n}=-K_\Sig D,\quad 0\le r''\le 2\
.\label{enn61}\end{equation}
\end{theorem}

This means that the number of real rational curves passing through
any generic conjugation-invariant configuration of points in
$\Sig$, where $r''$ is bounded as in the assertion, is
asymptotically equal to the number of all complex rational curves
in the logarithmic scale.

We propose a natural extension of Theorem \ref{tn4} to all Del
Pezzo surfaces and other values of $r''$:

\begin{conjecture}\label{con2} {\rm (\cite{IKS2})}
Let $D$ be an ample divisor on a real unnodal Del Pezzo surface
$\Sig$. Then, for any fixed $r''\ge 0$,
$$\lim_{n\to\infty}\frac{\log W_{r''}(\Sig,nD)}{n\log
n}=\lim_{n\to\infty}\frac{\log N_{\Sig,nD}}{n\log n}=-K_\Sig D\
.$$
\end{conjecture}

The following statement implies Conjecture \ref{con2} for the
hyperboloid:

\begin{theorem}\label{tnn10} {\rm (\cite{IKS2})} Let $D$ be an ample divisor on
$(\PP^1)^2$ of bi-degree $(a_1,a_2)$, $0<a_1\le a_2$, and let a
sequence of positive integers $r''(n)<(a_1+a_2)n$ satisfy
$\lim_{n\to\infty}r(n)/n=r''_0$. Then
\begin{equation}\lim_{n\to\infty}\inf\frac{\log
W_{r''(n)}((\PP^1)^2,nD)}{n\log
n}\ge-K_{(\PP^1)^2}D-r''_0=2a_1+2a_2-r''_0\
.\label{enn114}\end{equation}\end{theorem}

\medskip\noindent{\bf Organization of the material.}
The paper is structured as follows: in section \ref{secn13} we
introduce tropical curves, in section \ref{sec2} we describe
tropical limits of rational curves defined over the field $\K_\R$,
in section \ref{secn14} we prove the tropical formula for
Welschinger invariants, in section \ref{secn15} we obtain an
explicit combinatorial description for Welschinger invariants via
lattice paths and subdivisions of convex lattice polygons.

\medskip\noindent
{\bf Acknowledgements.} I am grateful to I. Itenberg and V.
Kharlamov for very valuable discussions. I thank the referee for a
careful reading of the manuscript and pointing out a number of
defects in the preliminary version.

\section{Non-Archimedean amoebas and tropical
curves}\label{secn13}

\subsection{Basic definitions and notation}\label{sec1} In
sections \ref{secn13} and \ref{secn14} we assume that $\Del$ is a
non-degenerate lattice polygon in $\R^2$, $\Sig=\Tor(\Del)$ is a
toric surface associated with $\Del$, and $|D|$ is the
tautological linear system, generated by the monomials $x^iy^j$,
$(i,j)\in\Del\cap\Z^2$.

\begin{figure}
\begin{center}
\epsfxsize 145mm \epsfbox{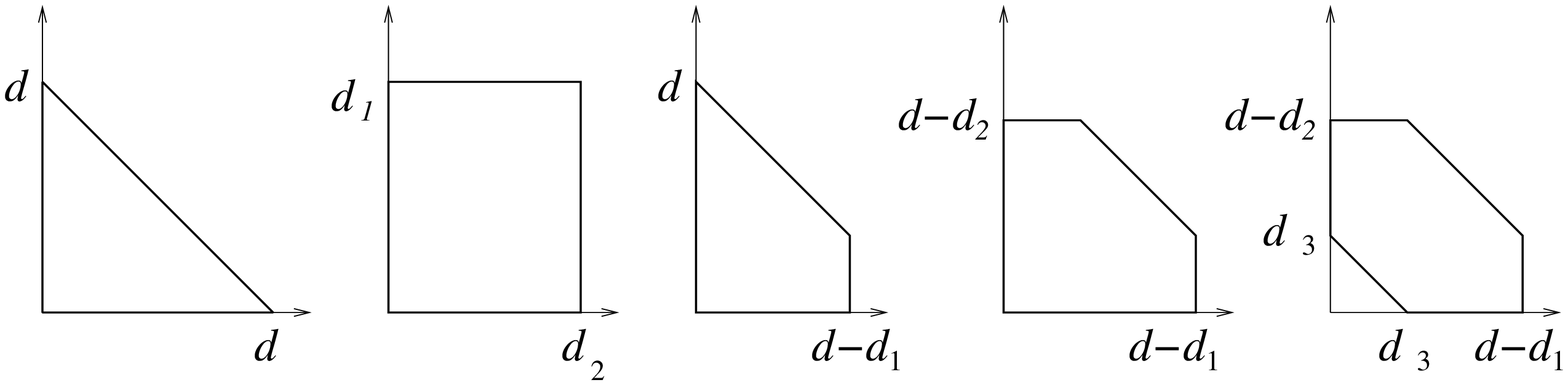}
\end{center}
\caption{Polygons associated with $\PP^2$, $(\PP^1)^2$, $\PP^2_1$,
$\PP^2_2$, and $\PP^2_3$} \label{fn2}
\end{figure}

In particular, the surfaces $\PP^2$, $(\PP^1)^2$, $\PP^2_k$,
$k=1,2,3$, are naturally associated with the polygons shown in
Figure \ref{fn2}:
\begin{itemize}\item the triangle $\conv\{(0,0),(0,d),(d,0)\}$, if
$\Sig=\PP^2$, $\deg D=d$,
\item the rectangle $\conv\{(0,0),(d_2,0),(d_2,d_1),(0,d_1)\}$, if
$\Sig=(\PP^1)^2$, $\deg D=(d_1,d_2)$, \item the trapeze
$\conv\{(0,0),(d-d_1,0),(d-d_1,d_1),(0,d)\}$ if $\Sig=\PP^2_1$,
$D\sim dL-d_1E_1$, \item the pentagon $\conv\{(0,0), (d-d_1,0),
(d-d_1,d_1),(d_2,d-d_2),(0,d-d_2)\}$, if $\Sig=\PP^2_2$, $D\sim
dL-d_1E_1-d_2E_2$, \item the hexagon $\conv\{(d_3,0), (d-d_1,0),
(d-d_1,d_1), (d_2,d-d_2), (0,d-d_2), (0,d_3)\}$ if $\Sig=\PP^2_3$,
$D\sim dL-d_1E_1-d_2E_2-d_3E_3$.
\end{itemize}

Observe that \begin{equation}r'+2r''=-K_\Sig
D-1=|\partial\Del\cap\Z^2|-1\ .\label{enn13}\end{equation}

The amoeba $A_C$ of a curve $C\in|D|_\K$ is defined as the closure
of the set $\val(C\cap(\K^*)^2)\subset\R^2$. By Kapranov's theorem
(see \cite{I,K,M}) the amoeba $A_C$ of a curve $C\in|D|_\K$, given
by the equation
\begin{equation}f(x,y):=\sum_{(i,j)\in\Del}a_{ij}x^iy^j=0\label{enn25}\end{equation} with $\Del$ as
the Newton polygon of $f$, is the corner locus of the convex
piece-wise linear function \begin{equation}N_f(x,y)=\max_{(i,j)\in
\Del\cap\Z^2}(xi+yj+\val(a_{ij})),\quad x,y\in\R\
.\label{enn1}\end{equation} In particular, $A_C$ is a planar graph
with all vertices of valency $\ge 3$, consisting of closed
segments and rays.

An amoeba $A$ with Newton polygon $\Del$ is called {\it
reducible}, if $A$ is the union of two amoebas $A',A''\ne A$ with
Newton polygons $\Del',\Del''$ such that $\Del=\Del'+\Del''$
(Minkowski sum).

Take the convex hull $\widetilde\Del$ of the set
$\{(i,j,-\val(a_{ij}))\in\R^3\ :\ (i,j)\in \Del\cap\Z^2\}$ and
define the function
\begin{equation}\nu_f:\Del\to\R,\quad\nu_f(x,y)=\min\{\gam\ :\
(x,y,\gam)\in\widetilde\Del\}\ .\label{enn4}\end{equation} This is
a convex piece-wise linear function, whose linearity domains form
a subdivision $S_C$\footnote{Clearly, $S_C$ does not depend on the
choice of a polynomial $f$ defining the curve $C$.} of $\Del$ into
convex lattice polygons $\Del_1,...,\Del_N$. The function $\nu_f$
is Legendre dual to $N_f$ (see, for instance, \cite{I}), and thus,
the subdivision $S_C$ is combinatorially dual to the pair
$(\R^2,A_C)$.

We define the {\it tropical curve}, corresponding to the algebraic
curve $C$, as a balanced graph, supported at $A_C$ (cf.
\cite{I,RST}), i.e., this is the non-Archimedean amoeba $A_C$,
whose edges are assigned weights equal to the lattice
lengths\footnote{We define the lattice length $|\sig|$ of a
segment $\sig$ with integral endpoints as $|\sig\cap\Z^2|-1$.} of
the dual edges of $S_C$. The subdivision $S_C$ can be uniquely
restored from the tropical curve $A_C$ \cite{RST} (we denote a
tropical curve and the supporting amoeba by the same symbol, no
confusion will arise).

\subsection{Rank of a tropical curve} Let $A$ be a
tropical curve with Newton polygon $\Del$, given by a tropical
polynomial (\ref{enn1}). On the right-hand side of formula
(\ref{enn1}) we remove unnecessary linear functions, take the
terms $\val(a_{ij})$ in the remaining linear functions as
variables, and factorizing by common shift, obtain the space
$\R^{|V(S)|-1}$, where $V(S)$ is the set of the vertices of the
dual subdivision $S$. The set $\Iso(A_C)$ of tropical curves with
Newton polygon $\Del$, which are combinatorially isotopic to the
given curve (or, equivalently, are dual to the same subdivision
$S$ of $\Del$), is defined in $\R^{|V(S)|-1}$ by certain linear
relations and inequalities, and thus forms a convex polyhedron.
Its dimension is called the {\it rank} of the tropical curve $A$.
Clearly, the rank of a tropical curve $\rk(A_C)$ is bounded from
below by the {\it virtual rank}
$$\rk_{vir}(A)=|\pi_0(\R^2\backslash A)|-1-\sum_{p\in
V(A)}(V_p-3)\ ,$$ where $V(A)$ is the set of vertices of $A$, and
$V_p$ is the valency of the vertex $p$. In terms of $S$,
$$\rk_{vir}(A)=\rk_{vir}(S)=|V(S)|-1-\sum_{k=1}^N(|V(\Del_k)|-3)\ ,$$ where
$V(S)$, $V(\Del_k)$ are the sets of vertices of $S$ and $\Del_k$,
respectively.

\begin{definition}\label{drev1} Let $\bx_1,...,\bx_m\in\Q^2$ be
distinct points, $A$ a tropical curve with Newton polygon $\Del$.
The set of tropical curves $A'\in\Iso(A)$ such that
$\bx_1,...,\bx_k$ are vertices of $A'$, and
$\bx_{k+1},...,\bx_m\in A'$ imposes a number of linear conditions
on the variables $\val(a_{ij})$, introduced above. We say, that
the pair of configurations
$(\{\bx_1,...,\bx_k\},\{\bx_{k+1},...,\bx_m\})$ is in $A$-generic
position if the aforemention set of curves $A'$ either is empty,
or has codimension $\ge 2k+(m-k)=m+k$ in $\Iso(A)$. We say that a
configuration $\bx_1,...,\bx_m$ is $\Del$-generic, if it is
generic for any division into a pair of configurations and any
tropical curve with Newton polygon $\Del$.
\end{definition}

\begin{lemma}\label{lnn4} The set of $\Del$-generic configurations of $m$
points in $\Q^2$ is dense in the set of all $m$-tuples in $\Q^2$.
In particular, if a tropical curve $A$ with Newton polygon $\Del$
passes through a $\Del$-generic configuration
$\bx_1,...,\bx_m\in\Q^2$ so that $\bx_1,...,\bx_k$ are vertices of
$A$, then
\begin{equation}m+k\le\rk(A)\ .\label{enn2}\end{equation}
\end{lemma}

\begin{proof} The requirements, imposed by the pair of configurations
$(\{\bx_1,...,\bx_k\},\{\bx_{k+1},...,\bx_m\})$ on $A'\in\Iso(A)$,
can be written as linear conditions on the variables
$\val(a_{ij})$ in $\R^{|V(S)|-1}$, two for a vertex $\bx_i$ of
$A'$, and one for $\bx_i\in A'$. The failure of generality in this
case means just a linear relation to the coordinates of
$\bx_1,...,\bx_m$. Since there are only finitely many
combinatorial types of tropical curves with Newton polygon $\Del$,
we obtain that the set of $\Del$-generic $m$-tuples is the
complement of finitely many hyperplanes in $(\Q^2)^m$. \end{proof}

A maximal straight line interval, contained in $A_C$, is called an
{\it extended edge} of $A_C$. The edges of $A_C$, forming an
extended edge, are dual to a sequence of parallel edges of $S_C$
ordered so that each two successive edges of this sequence belong
to one polygon $\Del_i$. We say that an extended edge of $A_C$ is
dual to each of the edges of $S_C$ in the corresponding sequence.

Having a tropical curve $A_C$ with Newton polygon $\Del$ and a
$\Del$-generic configuration of point on it, we call the set of
vertices of $A_C$, coinciding with points of the configuration,
and the set of extended edges of $A_C$, containing points of the
configuration in their interior, a {\it basic set of extended
edges and vertices} of $A_C$.

\section{Tropicalization of real rational algebraic
curves}\label{sec2} We work over the field $\K$ and its real
subfield $\K_\R=\bigcup_{m\ge 1}\R\{\{t^{1/m}\}\}$.

\subsection{Tropical limit} Fix some $\Del$-generic collections of points
$\overline\bx'=\{\bx'_1,...,\bx'_{r'}\}$ and
$\overline\bx''=\{\bx''_1,...,\bx''_{r''}\}$ in $\Q^2$. Let
$\bp'_1,...,\bp'_{r'}\in(\K_\R^*)^2$ be generic points, satisfying
$\val(\bp'_i)=\bx'_i$, $i=1,...,r'$, and let $\bp''_{ij}$,
$i=1,...,r''$, $j=1,2$, be generic points in
$(\K^*)^2\backslash(\K_\R^*)^2$, satisfying
$\conj(\bp''_{i1})=\bp''_{i2}$, $\val(\bp''_{ij})=\bx''_i$,
$i=1,...,r''$, $j=1,2$. That is
$$\bx'_i=(-\alp'_i,-\bet'_i),\
i=1,...,r',\quad\bx''_i=(-\alp''_i,-\bet''_i),\ i=1,...,r''\ ,$$
and \begin{equation}\bp'_i=(t^{\alp'_i}(\xi'_i+o(1)),\
t^{\bet'_i}(\eta'_i+o(1))),\quad \xi'_i,\eta'_i\in\R^*,\quad
i=1,...,r'\ ,\label{enn37}\end{equation}
\begin{equation}\begin{cases}\bp''_{i1}=(t^{\alp''_i}(\xi''_i+o(1)),\
t^{\bet''_i}(\eta''_i+o(1))),&\\
\bp''_{i2}=(t^{\alp''_i}(\overline{\xi''_i}+o(1)),\
t^{\bet''_i}(\overline{\eta''_i}+o(1))),&\end{cases}\quad\xi''_i,\eta''_i\in\C\backslash\R,
\quad i=1,...,r''\ .\label{enn38}\end{equation}

Let a rational curve $C\in|D|_\K$, given by a polynomial
\begin{equation}f(x,y):=\sum_{(i,j)\in\Del}a_{ij}(t)x^iy^j=0\ ,\label{enn3}\end{equation} be
defined over $\K_\R$, and pass through $\bp'_i$, $i=1,...r'$, and
$\bp''_{kj}$, $k=1,...,r''$, $j=1,2$.

Changing the parameter $t\mapsto t^m$, we make all the exponents
of $t$ in $a_{ij}(t)$, $(i,j)\in\Del$, integral and the function
$\nu_f$ integral-valued at integral points. Introduce the
polyhedron
$$\widetilde\Del=\{(\alp,\bet,\gam)\in\R^3\ :\
(\alp,\bet)\in\Del,\ \gam\ge\nu_f(\alp,\bet)\}\ .$$ It defines a
toric variety $Y=\Tor(\widetilde\Del)$, which is naturally fibred
over $\C$ so that the fibres $Y_t$, $t\ne 0$, are isomorphic to
$\Tor(\Del)$, and $Y_0$ is the union of toric surfaces
$\Tor(\widetilde\Del_i)$, $i=1,...,N$, with
$\widetilde\Del_1,...,\widetilde\Del_N$ being the faces of the
graph of $\nu_f$. By the choice of $\nu_f$,
$\Tor(\widetilde\Del_i)\simeq\Tor(\Del_i)$, and we shall write
that $Y_0=\bigcup_i\Tor(\Del_i)$. Then (\ref{enn3}) defines an
analytic surface $\{f(x,y)=0\}$ in a neighborhood of $Y_0$, which,
by \cite{ShP}, Lemma 2.3, fibers into equisingular rational curves
$C^{(t)}\subset Y_t\simeq\Tor(\Del)$, and whose closure intersects
$Y_0$ along the curve $C^{(0)}$ that can be identified with
$\bigcup_iC_i\subset\bigcup_i\Tor(\Del_i)$. Passing if necessary
to a finite cyclic covering ramified along $Y_0$, we can make
$\Tor(\widetilde\Del)$ non-singular everywhere but may be at
finitely many points, corresponding to the vertices of
$\widetilde\Del$, and, in addition, make the surfaces
$\Tor(\Del_k)\backslash\Sing(\Tor(\widetilde\Del))$, $k=1,...,N$,
smooth and intersecting transversally in
$\Tor(\widetilde\Del)\backslash\Sing(\Tor(\widetilde\Del))$. As in
\cite{ShP}, we define the tropical limit (tropicalization) of the
curve $C$ to be the pair, consisting of the tropical curve $A_C$,
and a collection of real curves $C_k=C^{(0)}\cap\Tor(\Del_k)$,
$k=1,...,N$, which are defined by
\begin{equation}f_k(x,y):=\sum_{(i,j)\in\Del_k}a_{ij}^0x^iy^j=0,\quad
k=1,...,N\ ,\label{enn7}\end{equation} respectively, where
$a_{ij}(t)=(a_{ij}^0+O(t))t^{\nu_f(i,j)}$, $(i,j)\in\Del$.

\subsection{Tropical limits of real rational curves}\label{secn2}
A tropical curve $A_C$ is called {\it nodal} if the polygons of
the dual subdivision $S_C$ are triangles and parallelograms.
Recall that, by \cite{ShP}, Lemma 2.2, $\rk(A_C)=\rk_{vir}(A_C)$
if $A_C$ is nodal. Notice also that the weights of the edges of a
nodal tropical curve $A_C$ are constant along its extended edges,
and thus, one can speak of weights of extended edges of nodal
tropical curves.

An irreducible nodal tropical curve $A_C$ is called {\it real
rational of type} $(r',r'',s'')$ with $r',r''$, satisfying
(\ref{enn36}), and $0\le s''\le r''$, if
\begin{itemize}
\item
$\rk(A_C)=\rk_{vir}(A_C)=|\partial\Del\cap\Z^2|-1-r''+s''$,
\item the weights of the semi-infinite edges of
$A_C$ are $1$ or $2$, i.e., the edges of $S_C$ lying on
$\partial\Del$ are of length $1$ or $2$.
\end{itemize}

\begin{proposition}\label{pn1}
Under the assumptions of section \ref{sec2}, $A_C$ is a real
rational tropical curve of type $(r',r'',s'')$ for some $0\le
s''\le r''$, passing through $\bx'_1,...,\bx'_{r'}$ and
$\bx''_1,...,\bx''_{r''}$ in such a way that precisely $s''$ of
the points $\bx''_1,...,\bx''_{r''}$ are trivalent vertices of
$A_C$, whereas $\bx'_1,...,\bx'_{r'}$ and the remaining $r''-s''$
points among $\bx''_1,...,\bx''_{r''}$ (which we denote
$\bw_1,...,\bw_{r''-s''}$) are not vertices of $A_C$, and,
moreover, $\bw_1,...,\bw_{r''-s''}$ lie on edges of $A_C$ of even
weight. Furthermore, for $k=1,...,N$,
\begin{itemize}\item if $\Del_k$ is a parallelogram, then $f_k$ is the product of a monomial and few irreducible
binomials, \item if $\Del_k$ is a triangle, then $C_k$ is either a
real rational curve crossing
\mbox{$\Tor(\partial\Del_k):=\bigcup_{\sig\subset\partial\Del_k}\Tor(\sig)$}
at $3$, $4$, or $5$ points, at which it is unibranch, or $C_k$ is
the union of two imaginary conjugate rational curves such that any
component of $C_k$ crosses $\Tor(\partial\Del_k)$ at precisely
three points and is unibranch there.
\end{itemize}
\end{proposition}

\begin{proof} Our argumentation is similar to that in \cite{ShP}, section
3.3, used to establish that the non-Archimedean amoebas of nodal
curves in toric surfaces passing through a respective number of
generic points, are nodal.

\medskip
{\it Step 1}. Let $\bx'_1,...,\bx'_{s'}$ and
$\bx''_1,...,\bx''_{s''}$ be vertices of $A_C$ for some $0\le
s'\le r'$, $0\le s''\le r''$, and $\bx'_i$, $i>s'$, $\bx''_j$,
$j>s''$, not be vertices. By (\ref{enn2}),
\begin{equation}r'+r''+s'+s''\le\rk(A_C)\ .\label{enn17}\end{equation}

By \cite{ShP}, Lemma 2.2, we have
$$\rk(A_C)=\rk_{vir}(A_C)+d(A_C)=|V(S_C)|-1-\sum_{k=1}^N(|V(\Del_k)|-3)
+d(A_C)\ ,$$ where
$$2d(A_C)\le\begin{cases}0,\quad \text{if}\
A_C\ \text{is nodal},&\\
\sum_{m\ge 2}((2m-3)N_{2m}-N'_{2m})+\sum_{m\ge
2}(2m-2)N_{2m+1}-1,\quad \text{otherwise},&\end{cases}$$ $N_l$
being the number of $l$-gons in $S_C$, and $N'_{2m}$ being the
number of $2m$-gons in $S_C$ whose opposite edges are parallel.
Substituting the two latter relations into (\ref{enn17}) and using
the Euler formula for the subdivision $S_C$, one obtains
\begin{equation}2(r'+r''+s'+s'')\le|V(S_C)\cap\partial\Del|+\begin{cases}N_3,\quad\text{if}\
A_C\ \text{is nodal},&\\ N-\sum_{m\ge
2}N'_{2m}-1,\quad\text{otherwise}.&\end{cases}\label{enn19}\end{equation}

Denote by $\hat\chi$ the Euler characteristic of the
normalization. Let $C_{ij}$, $1\le j\le n_i$, be all the
components of $C_i$, $1\le i\le N$, repeating each component
$C_{ij}$ with its multiplicity in $C_i$, and let $s_{ij}$ be the
number of local branches of $C_{ij}$ centered along
$\Tor(\partial\Del_i)$. Since the rational curves $C^{(t)}$, $t\ne
0$, degenerate into $C^{(0)}$ in a flat family, all the components
$C_{ij}$ of $C^{(0)}$ are rational in view of the inequality for
geometric genera $g(C^{(t)})\ge\sum_{i,j}g(C_{ij})$ (see
\cite{DH}, Proposition 2.4, or \cite{N}\footnote{The hypotheses of
Proposition 2.4 \cite{DH} require that all members of the family
are reduced curves, which can be achieved in our situation by the
normalization of the whole family $C$.}). Denote by ${\mathcal
C}_b$ the set of components $C_{ij}$, defined by irreducible
binomials, and by ${\mathcal C}_{nb}$ the set of remaining
components $C_{ij}$. Let $U$ be the union of regular neighborhoods
in the three-fold $Y$ of the intersection points of $C^{(0)}$ with
$\bigcup_k\Tor(\partial\Del_k)$. Then
$$2=\hat\chi(C^{(t)})\le\sum_{k=1}^N\sum_{C_{kj}\in{\mathcal C}_{nb}}(2-s_{kj})+\hat\chi(C^{(0)}\cap
U)$$\begin{equation} \le\sum_{k=1}^N\sum_{C_{kj}\in{\mathcal
C}_{nb}}(2-s_{kj})+s(\partial\Del)\ ,\label{enn18}\end{equation}
the latter inequality following from \cite{ShP}, Lemma 3.2 and
Remark 3.4, where $s(\partial\Del)$ stands for the number of local
branches (counting multiplicities) of the curves $C_1,...,C_N$
centered at $\Tor(\sig)$ with $\sig$ running over all the edges of
$\Del$.

\medskip
{\it Step 2}. For an estimation of the right-hand side of
(\ref{enn18}) we shall construct some graph $G$.

First, we compactify $\R^2$ into $S^2$ by adding an infinite point
$v_\infty$. The tropical curve $A_C$ is then compactified by
closing the semi-infinite edges with the point $v_\infty$.

For any $k=1,...,N$, the reduced curve $C_k^\red$ is split into
irreducible real components and pairs of imaginary conjugate
irreducible components. All such real components or pairs of
imaginary conjugate components for all $k=1,...,N$, except for
real irreducible components from the set ${\mathcal C}_b$, and the
point $v_\infty$ are taken as the vertices of the graph
$\widetilde G$. If $\Del_k\cap\Del_l=\sig$ is a common edge, and
components $C'_k$ or $C_k$ and $C'_l$ of $C_l$ contain a common
pair of imaginary conjugate points in $\Tor(\sig)$, then we join
$C'_k$ and $C'_l$ by an edge. If $\sig\subset\partial\Del$ is an
edge of $\Del_k$, and a component $C'_k$ of $C_k$ contains a pair
of imaginary conjugate points in $\Tor(\sig)$, then we join $C'_k$
and $v_\infty$ by an edge.

The constructed graph $\widetilde G$ is then transformed as
follows. Let $C'_k$ be a part of $C_k$, consisting of two
imaginary conjugate components from the set ${\mathcal C}_b$. The
curve $C'_k$ intersects with $\Tor(\sig_1)$ and $\Tor(\sig_2)$,
where $\sig_1,\sig_2$ is a pair of parallel edges of $\Del_k$. If
$\Del_k\cap\Del_p=\sig_1$, $\Del_k\cap\Del_q=\sig_2$, then we
remove the vertex $C'_k$ and all ending at it edges, and instead
join by an edge any two curves $C'_p$ and $C'_q$ which have been
joined with $C'_k$. If $\Del_p\cap\Del_k=\sig_1$ and
$\sig_2\subset\partial\Del$, then we get rid of the vertex $C'_k$
of the graph and all adjacent edges, instead connecting $v_\infty$
with any curve $C'_p$, which has been joined with $C'_k$. In this
manner, we get rid of one-by-one all the vertices of $\widetilde
G$, corresponding to the curves $C'_k$ which consist of two
imaginary conjugate components from the set ${\mathcal C}_b$.

We observe that the graph $\widetilde G$ naturally projects to
$A_C\cup\{v_\infty\}$, when sending any vertex $C'_k$ of
$\widetilde G$ to the vertex $v$ of $A_C$ dual to $\Del_k$, and
sending edges of $\widetilde G$ to respective segments of the
compactified extended edges of $A_C\cup\{v_\infty\}$.

The required graph $G$ will be a subgraph of $\widetilde G$.

First, we define a subgraph $G_0$ of $\widetilde G$ as follows.
Each point $\bx''_j$, $j>s''$, lies on an edge of $A_C$. Let
$\sig$ be the dual edge of $S_C$. By \cite{ShP}, formula (3.7.17),
$f_k^\sig(\xi''_j,\eta''_j)=f_k^\sig(\overline
\xi''_j,\overline{\eta''_j})=0$, where $f_k^\sig$ is the
truncation of the polynomial $f_k(x,y)$, defined by (\ref{enn7}),
on the edge $\sig$, and $\xi''_j,\eta''_j$ are taken from formula
(\ref{enn38}). That is, $\Tor(\sig)$ contains a pair of distinct
imaginary conjugate points, $z$ and $\conj(z)$, which belong to
curves $C_k$ such that $\Del_k\supset\sig$. Hence the edge of
$A_C$, containing $\bx''_j$, is covered by the projection of some
edges of $\widetilde G$. We choose one such covering edge of
$\widetilde G$, denote it by $\eps_j$, and build the graph $G_0$
from the edges $\eps_j$, $j>s''$, with their endpoints as vertices
of $G_0$.

Notice that due to the $\Del$-generic position of the points
$\bx''_j$, $j>s''$, the graph $G_0\backslash\{v_\infty\}$ is the
union of trees, the valency of its vertices that differ from
$v_\infty$ is at most two, and the intersection of the projections
of any two of its edges is finite.

A vertex $v\ne v_\infty$ of $G_0$ of valency two corresponds to
some curve $C'_k$, $1\le k\le N$. If $C'_k$ is the union of two
imaginary conjugate rational curves such that each of them crosses
$\Tor(\partial\Del_k)$ at precisely three points, and is unibranch
at these intersection points, then we call the vertex $v$ of $G_0$
an {\it extendable vertex}. Notice that the components of $C'_k$
are defined by polynomials with the same Newton triangle, and are
determined by $\bx''_j$, $j>s''$, up to a finite choice in view of
\cite{ShP}, Lemma 3.5. Furthermore, due to the generality of the
coefficients $\xi''_i,\eta''_i$ in (\ref{enn38}), the intersection
$C'_k\cap\Tor(\sig)$ is either empty or a pair of imaginary
conjugate points for any $\sig\subset\partial\Del_k$. Exactly one
of these pairs lies on an edge of $\Del_k$, whose dual $A_C$-edge
is not covered by $G_0$. However, the latter $A_C$-edge is covered
by some edges of $\widetilde G$ with an endpoint $C'_k$. We choose
one such edge of $\widetilde G$ and append it to $G_0$.

Performing this procedure for all extendable vertices of $G_0$, we
obtain a graph $G_1$. Next we determine the extendable vertices of
$G_1$, using the same definition, and append, if necessary, new
edges to $G_1$, obtaining $G_2$. Repeating this procedure, we
finally end up with some subgraph $G$ of $\widetilde G$.

\medskip
{\it Step 3}. An important observation is that no new edge in
$G_1\backslash G_0$ can join two extendable vertices of $G_0$.
Indeed, the position of two non-adjacent extendable vertices of
$G_0$ is uniquely determined by the position of the corresponding
four points among $\bx''_1,...,\bx''_{s''}$, and thus, due to the
generality of the latter points, the straight line through the
given extendable vertices is not orthogonal to any of the segments
joining integral points in $\Del$. By a similar reason, no edge
from $G_{i+1}\backslash G_i$ can join two extendable vertices of
$G_i$, $i\ge 1$. In particular, this means that the number of
edges in $G\backslash G_0$ is equal to the total number of
extendable vertices in $G_0,G_1,...$.

Furthermore, assume that $C'_k$ is a $q$-valent vertex of $G$,
$q\ge 3$, which is not extendable for any $G_i$, $i\ge 0$. Then
$C'_k$ cannot be the union of two imaginary conjugate rational
curves, such that each of them crosses $\Tor(\partial\Del_k)$ at
precisely $q$ points and is unibranch there. Indeed, if the
situation were as described, the position of any $q-1$ of the
intersection points of a component of $C'_k$ with $\Tor(\Del_k)$
would determine this component and the fourth intersection point
uniquely up to a finite choice. Thus, the relation on the
coordinates of the intersection points of $C'_k$ with
$\Tor(\partial\Del_k)$ would imply a relation on the coordinates
of the points $\bp''_{j1}$, $j>s''$, contrary to their general
choice. In particular, this implies that the contribution of such
a vertex to the sum on the right-hand side of (\ref{enn18}) is
$\le2-2q$, and the equality corresponds to the case of $C'_k$
being the union of two imaginary conjugate components such that a
component crosses $\Tor(\partial\Del_k)$ at $q+1$ points and is
unibranch there.

We also observe that the points $\bx'_1,...,\bx'_{s'}$ and
$\bx''_1,...,\bx''_{s''}$, which are vertices of $A_C$ by
assumption, are not projections of the vertices of $G$. Indeed,
the points $\bx'_1,...,\bx'_{s'}$ and $\bx''_1,...,\bx''_{s''}$
cannot be projections of the vertices of $G_0$ in view of a
general position of these points with respect to $\bx''_j$,
$j>s''$, which in turn determine the projections of the edges of
$G_0$. Next, the projections of the edges of $G_1$, and thus, the
projections of the extendable vertices of $G_1$ are determined by
$G_0$ up to a finite choice, where the number of choices is
bounded from above by the total number of possible lattice edges
in $\Del$. By induction we get that the projections of the edges
and vertices of $G$ are determined by $\bx''_j$, $j>s''$, up to a
finite choice, and hence the general position of the points
$\bx'_1,...,\bx'_{s'}$ and $\bx''_1,...,\bx''_{s''}$ does not
allow them to be the projections of the vertices of $G$.

\medskip
{\it Step 4}. Denote by $N_q(G)$, $q\ge 1$, the number of
$q$-valent vertices of $G$ different from $v_\infty$. Denote by
$N_3^e(G)$ the number of trivalent vertices of $G$ which are
extendable for some $G_i$, $i\ge 0$. At last, denote by $V_\infty$
the valency of the vertex $v_\infty$ of $G$. Then the contribution
of a vertex of $G$ to the right-hand side of (\ref{enn18}) is
\begin{itemize}\item $\le -q-1$, if it is a vertex of valency $q=1,2$,
\item $\le 2-2q$, if it is a vertex of valency $q\ge 3$, which is not extendable for any $G_i$, $i\ge 0$, \item $-2$,
if it is a trivalent vertex, which is extendable for some $G_j$,
$j\ge 0$,
\end{itemize} and thus, the total contribution of the vertices of
$G$ is $$\le-\sum_{q\ge 1}(q+1)N_q(G)-\sum_{q\ge
4}(2q-2)N_q(G)+2N_3^e(G)\ .$$ Since the total number of edges of
$G$ is $r''-s''+N_3^e(G)$, we can rewrite the latter inequality as
$$\le-\sum_{q\ge 1}N_q(G)-\sum_{q\ge
4}(2q-2)N_q(G)-2r''+3s''+V_\infty\ .$$

Denote by $\widetilde N$ the number of polygons in $S_C$, which
are dual to the vertices of $A_C$ not covered by the vertices of
$G$, and which are not even-gons, whose pairs of opposite sides
are parallel. Clearly, a curve $C_k$ corresponding to such a
polygon $\Del_k$, contains a component from the set ${\mathcal
C}_{nb}$, and thus, contributes $\le -1$ to the sum on the
right-hand side of (\ref{enn18}). That is, inequality
(\ref{enn18}) implies
$$2\le-\sum_{q\ge 1}N_q(G)-\sum_{q\ge
4}(2q-2)N_q(G)-2r''+3s''+V_\infty-\widetilde N+s(\partial\Del)$$
\begin{equation}\Longrightarrow\quad2\le-N-\sum_{q\ge
4}(2q-2)N_q(G)+\sum_{m\ge
2}N'_{2m}-2r''+2s''+V_\infty+s(\partial\Del)\
,\label{enn21}\end{equation} with an equality only if
\begin{itemize}\item the vertices of $G$ injectively project to the set of
vertices of $A_C$, which are dual to polygons different from
even-gons, whose all pairs of opposite sides are parallel; \item
if $\Del_k$ is an odd-gon, whose dual vertex of $A_C$ is not
covered by the vertices of $G$, then $C_k$ contains precisely one
(counting multiplicities) component from ${\mathcal C}_{nb}$, and,
moreover, this component crosses $\Tor(\partial\Del_k)$ at
precisely three points and is unibranch there.\end{itemize}
Combining (\ref{enn19}) and (\ref{enn21}), we derive
$$2r'+4r''+2s'+\sum_{q\ge
4}(2q-2)N_q(G)\le-2+|V(S_C)\cap\partial\Del|+s(\partial\Del)+V_\infty+\begin{cases}
0,\ \text{if}\ A_C\ \text{is nodal},&\\ -1,\
\text{otherwise}&\end{cases}$$
$$\stackrel{\text{(\ref{enn13})}}{\Longrightarrow}\quad\sum_{q\ge
4}(2q-2)N_q(G)2s'+|\partial\Del\cap\Z^2\backslash
V(S_C)|+(|\partial\Del\cap\Z^2|-s(\partial\Del))$$\begin{equation}\le
V_\infty+\begin{cases} 0,\ \text{if}\ A_C\ \text{is nodal},&\\
-1,\ \text{otherwise}.&\end{cases}\label{enn22}\end{equation} The
immediate relation
\begin{equation}V_\infty\le|\partial\Del\cap\Z^2\backslash V(S_C)|\label{enn44}\end{equation} excludes
the case of a non-nodal tropical curve $A_C$, and, in the case of
a nodal tropical $A_C$, implies the equality in (\ref{enn22}) and
in (\ref{enn44}) with the restriction $q\le 3$ on the valency of
vertices of $G$.

We then collect all the equality conditions mentioned above into
the following restrictions: \begin{enumerate}\item[(R1)] relations
(\ref{enn17}) and (\ref{enn18}) are equalities with $s'=0$ and
$s(\partial\Del)=|\partial\Del\cap\Z^2|$;\item[(R2)] the edges of
$S_C$ lying in $\partial\Del$ have length $1$ or $2$;\item[(R3)]
the semi-infinite extended edges of $A$, which are dual to the
edges of $S$ of length $>1$ lying on $\partial\Del$, belong to the
graph $G$;
\item[(R4)] any curve $C_k$, $1\le
k\le N$, crosses the divisors $\Tor(\sig)$,
$\sig\subset\partial\Del_k\cap\partial\Del$ only at its
non-singular points, and all these intersections are transversal;
\item[(R5)] all the components of the curves $C_k$ with Newton
parallelograms $\Del_k$ belong to ${\mathcal C}_b$; \item[(R6)]
the curves $C_k$, $1\le k\le N$, such that $\Del_k$ is a triangle
dual to a vertex of $A_C$ outside $G$, are rational, intersecting
$\Tor(\partial\Del_k)$ at precisely three points being unibranch
there; \item[(R7)] if $\Del_k$, $1\le k\le N$, is dual to a
univalent vertex of $G$, then $\Del_k$ is a triangle, $C_k$ is
rational and crosses $\Tor(\partial\Del_k)$ at precisely four
points being unibranch there; \item[(R8)] if $\Del_k$, $1\le k\le
N$, is dual to a bivalent vertex of $G$, then $\Del_k$ is a
triangle, $C_k$ is rational and crosses $\Tor(\partial\Del_k)$ at
precisely five points being unibranch there; \item[(R9)] if
$\Del_k$, $1\le k\le N$, is dual to a trivalent vertex of $G$,
extendable for some $G_i$, then $\Del_k$ is a triangle, $C_k$ is
the union of two imaginary conjugate rational curves, each of them
crossing $\Tor(\partial\Del_k)$ at precisely three points being
unibranch there.
\end{enumerate}

\medskip
{\it Step 5}. Now we show that $A_C$ is irreducible, and that any
trivalent vertex of $G$ is extendable for some $G_i$, $i\ge 0$.

First, the equality in (\ref{enn18}) means that, in the
deformation $C^{(t)}$, $t\ge 0$, no intersection point of any two
distinct components of $C_i$, $1\le i\le N$, is smoothed out.
Assuming that $A_C$ is reducible, and using the description of the
curves $C_k$, obtained in Step 4, we immediately conclude that
then the curve $C^{(t)}$, $t\ne 0$, becomes reducible.

Second, the existence of a trivalent vertex of $G$, which was not
extendable for all $G_i$, $i\ge 0$, would mean a non-trivial
relation to the position of the points $\bx''_j$, $j>s''$,
contrary to their $\Del$-generality.

\medskip
{\it Step 6}. To complete the proof of Proposition \ref{pn1}, it
remains to verify that the points $\bx''_1,...,\bx''_{s''}$ cannot
be four-valent vertices of $A_C$.

Assume on the contrary that, say, $\bx''_1$ is a four-valent
vertex of $A_C$, which is dual to a parallelogram $\Del_k$. The
points $\bp''_{11}$ and $\bp''_{12}$ pick out two conjugate
binomial components of the curve $C_k$, which cross two toric
divisors of $\Tor(\Del_k)$, corresponding to a pair of parallel
sides of $\Del_k$. Let these sides be dual to an extended edge
$\sig'$ of $A_C$; denote by $\sig''$ the other extended edge of
$A_C$, passing through $\bx''_1$. We claim that $A_C$ can be
deformed inside $\Iso(A_C)$ so that the deformed tropical curves
will pass through the configuration
$\overline\bx'\cup\overline\bx''$, the points
$\bx''_2,...,\bx''_{s''}$ will remain their vertices, and the edge
$\sig''$ will move out of $\bx''_1$. Indeed, by (\ref{enn17}),
reducing the requirement that $\sig''$ passes through $\bx''_1$,
we obtain at least one-dimensional family of curves in
$\Iso(A_C)$, whereas, keeping that requirement we obtain only the
curve $A_C$. By Lemma \ref{lnn11}(1), presented below in section
\ref{secn14}, any of the above deformed tropical curves equipped
with the same limit curves $C_1,...,C_N$ gives rise to a real
rational curve in $|D|_\K$, passing through $\bp'_i$, $i=1,...r'$,
and $\bp''_{kj}$, $k=1,...,r''$, $j=1,2$. Thus, we obtain
infinitely many such curves in $|D|_\K$, contradicting to the
generality of $\bp'_i$, $i=1,...r'$, and $\bp''_{kj}$,
$k=1,...,r''$, $j=1,2$. \end{proof}

\begin{remark}\label{rnn1}
We point out that, given a real tropical curve as in the assertion
of Proposition \ref{pn1}, and a collection of curves
$C_1,...,C_N$, as described in the proof, the curves $C^{(t)}$,
$t\ne 0$, which may appear from these data, are irreducible
provided that the graph $G$ differs from $A_C$. The latter holds
if, for example, $r'+s''>0$.
\end{remark}

Given a real rational tropical curve $A$ of type $(r',r'',s'')$
with Newton polygon $\Del$, passing through the points
$\bx'_1,...,\bx'_{r'}$ and $\bx''_1,...,\bx''_{r''}$ as described
in Proposition \ref{pn1}, to decide what are the curves
$C_1,...,C_N$, we have to know the graph $G$, which appeared in
the proof of Proposition \ref{pn1}.

We can reconstruct $G$ in the following way. Take all the extended
edges of $A$ passing through those points of
$\bx''_1,...,\bx''_{r''}$, which are not vertices of $A$. Recall
that all these edges are of even weight. Their union together with
endpoints (including $v_\infty$ if necessary) constitutes the
graph $G_0$. Let $V^2(G_0)$ be the set of bivalent vertices of
$G_0$, different from $v_\infty$. Take any subset $V_0\subset
V^2(G_0)$ and, for each $v\in V_0$, add to $G_0$ the extended edge
of $A$ with endpoint $v$, which is not in $G_0$. The new graph is
denoted by $G_1$. Then we take any subset $V_1\subset
V^2(G_1)\backslash V^2(G_0)$, and append the extended edges to
$G_1$ in a similar manner. In finitely many steps we end up with
some graph $G$, whose vertices are some trivalent vertices of $A$
and, perhaps, $v_\infty$, whose edges are some extendable edges of
$A$. Notice that all the edges of $G$ necessarily have even
weight. By the restrictions obtained in Step 4 of the proof of
Proposition \ref{pn1}, the set of trivalent vertices of $G$ must
be $\bigcup_{i\ge 0}V_i$.

\subsection{Real rational algebraic curves in the tropical limits of rational curves over $\K_\R$}

\begin{lemma}\label{lnn5} In the notation of section \ref{secn2},
let $\Del_k$ be a triangle in $S_C$, dual to a vertex of $A_C$,
which differs from $\bx''_1,...,\bx''_{r''}$, and is not a vertex
of $G$. Then $C_k$ is a real rational nodal curve, crossing
$\Tor(\partial\Del_k)$ at precisely three points, where it is
non-singular. The number of real solitary nodes of $C_k$ has the
same parity as $|\Int(\Del_k)\cap\Z^2|$. Furthermore, if the
lengths of edges of $\Del_k$ are odd, and we are given the
coefficients of the polynomial $f_k$ at the vertices of $\Del_k$,
then $f_k$ and thereby $C_k$ are determined uniquely.
\end{lemma}

\begin{proof} The statements immediately follow from \cite{ShP}, Lemma 3.5
and Proposition 6.1. \end{proof}

\begin{lemma}\label{lnn6}
In the notation of section \ref{secn2}, let $\Del_k$ be a triangle
in $S_C$, dual to a vertex $\bx''_j$, $j\le s''$, of $A_C$. Then
$C_k$ is a real rational nodal curve, crossing
$\Tor(\partial\Del_k)$ at precisely three points, where it is
non-singular. The number of real solitary nodes of $C_k$ has the
same parity as $|\Int(\Del_k)\cap\Z^2|$. Furthermore, if the
lengths of edges of $\Del_k$ are odd, then there are precisely
$|\Del_k|$ real curves $C_k$ as above, where $|\Del_k|$ is the
double Euclidean area of $\Del_k$.
\end{lemma}

\begin{proof} The fact that $C_k$ is non-singular along
$\Tor(\partial\Del_k)$, is nodal, and the number of solitary real
nodes is $|\Int(\Del_k)\cap\Z^2|$ $\mod 2$, follows from
\cite{ShP}, Lemma 3.5 and Proposition 6.1. To establish the last
statement of the lemma, we perform an affine transformation of the
lattice $\Z^2$, which takes $\Del_k$ into the triangle with
vertices $(0,m_1)$, $(m_2,0)$, $(m_3,0)$ with odd $m_1,m_2$ and
even $m_3$ such that $0\le m_2<m_3$, $0<m_1<m_3$. Then $C_k$ has
parameterization $x(\rho)=a\rho^{m_1}$,
$y(\rho)=b\rho^{m_2}(\rho-1)^{m_3-m_2}$ in suitable affine
coordinates $x,y$, where $a,b\in\R\backslash\{0\}$. The condition
that $\Del_k$ is dual to $\bx''_j$ means that there is $\tau\in\C$
such that $x(\tau)=\xi''_j$, $y(\tau)=\eta''_j$, where
$\xi''_j,\eta''_j$ are defined by (\ref{enn38}). Due to the
generality of $\xi''_j,\eta''_j$, the two latter equations have
$m_1(m_3-m_2)=|\Del_k|$ distinct solutions
$(a,b,\tau)\in\R^2\times\C$. \end{proof}

\begin{lemma}\label{lnn7}
In the notation of section \ref{secn2}, let $\Del_k$ be a triangle
in $S_C$, dual to a univalent vertex of the graph $G$. Then $C_k$
is a real rational nodal curve, crossing $\Tor(\partial\Del_k)$ at
precisely four points, two real and two imaginary conjugate, where
it is non-singular. The number of real solitary nodes of $C_k$ is
equal to $|\Int(\Del_k)\cap\Z^2|$ $\mod 2$. Assume that $\Del_k$
has two edges of odd length and one edge of even length, and we
are given the intersection points of $C_k$ with
$\Tor(\partial\Del_k)$. Then there are $|\Del_k|/|\sig|$ such
curves $C_k$, where $\sig\subset\Del_k$ is the edge of even
length.
\end{lemma}

\begin{proof} A suitable affine automorphism of the lattice $\Z^2$ takes
$\Del_k$ into the triangle with vertices $(0,m_1)$, $(2m_2,0)$,
$(2m_3,0)$, where $0\le m_2<m_3$, $m_1<2m_3$. We choose a
parameterization $x(\rho),y(\rho)$ of $C_k$ such that the value
$\rho=\infty$ corresponds to the intersection point of $C_k$ with
$\Tor([(0,m_1),(2m_3,0)])$, and the values $\rho=\pm\sqrt{-1}$
correspond to the intersection points with
$\Tor([(2m_2,0),(2m_3,0)])$. Then
$$x(\rho)=a(\rho-c)^{m_1},\quad
y(\rho)=b(\rho-c)^{2m_2}(\rho^2+1)^{m_3-m_2}\ ,$$ where
$a,b,c\in\R^*$. This immediately proves that $C_k$ is non-singular
along $\Tor(\partial\Del_k)$, as well as the fact that $C_k$ is
nodal and has no real non-solitary nodes (cf. \cite{ShP},
Proposition 6.1). Assume now that $m_1$ is odd. The conditions of
the fixed intersection points with $\Tor(\partial\Del_k)$ reduce
to the system of equations
$$\frac{a}{b(c^2+1)^{m_3-m_2}}=\const\in\R^*,\quad a(\sqrt{-1}-c)^{m_1}=\const\in\C^*\ ,$$
which has $m_1=|\Del_k|/|\sig|$ distinct solutions
$(a,b,c)\in(\R^*)^3$. \end{proof}

\begin{lemma}\label{lnn8}
In the notation of section \ref{secn2}, let $\Del_k$ be a triangle
in $S_C$, dual to a bivalent vertex of the graph $G$. Then $C_k$
is a real rational nodal curve, crossing $\Tor(\partial\Del_k)$ at
precisely five points, one real and two pairs of imaginary
conjugate points, where it is non-singular.
\end{lemma}

\begin{proof} As in the proof of Lemma \ref{lnn7}, we can assume that
$\Del_k=\conv\{(0,2m_1),(2m_2,0),(2m_3,0)\}$ with integers $0\le
m_2<m_3$, $m_1<2m_3$. We choose a parameterization
$x(\rho),y(\rho)$ of $C_k$ so that the value $\rho=\infty$
corresponds to the intersection point of $C_k$ with
$\Tor([(0,m_1),(2m_3,0)])$, the values $\rho=\pm\sqrt{-1}$
correspond to the intersection points with
$\Tor([(2m_2,0),(2m_3,0)])$, and the values
$\rho=\kappa,\overline\kappa$ for some
$\kappa\in\C\backslash(\R\cup\{\pm\sqrt{-1}\})$ correspond to the
intersection points with $\Tor([(0,2m_1),(2m_2,0)])$. Then
$$x(\rho)=a(\rho-\kappa)^{m_1}(\rho-\overline\kappa)^{m_1},\quad
y(\rho)=b(\rho-\kappa)^{m_2}(\rho-\overline\kappa)^{m_2}(\rho^2+1)^{m_3-m_2}\
,$$ with some $a,b\in\R^*$. It is immediate that $C_k$ is
non-singular at $C_k\cap\Tor(\partial\Del_k)$. Along the proof of
Proposition \ref{pn1}, the two pairs of imaginary intersection
points of $C_k$ with $\Tor(\partial\Del_k)$ are generic, which
reads as a system of equations on $a,b,\kappa$:
$$\frac{b}{a}(\kappa-\overline\kappa)^{m_2-m_1}(\kappa^2+1)^{m_3-m_2}=\lam,\quad
a(\sqrt{-1}-\kappa)^{m_1}(\sqrt{-1}-\overline\kappa)^{m_1}=\mu$$
with some generic $\lam,\mu\in\C\backslash\R$. This results in a
finitely many generic values for $\kappa$. In particular, we
obtain that $x'(\rho),y'(\rho)$ do not vanish simultaneously in
$(\C^*)^2\subset\Tor(\Del_k)$, or, equivalently, $C_k$ is an
immersed line $\C P^1$ into $\Tor(\Del_k)$. The condition that
three local branches of $C_k$ are centered at some point of
$(\C^*)^2$ reduces to a system of four equations
$$\begin{cases}&(\rho_1-\kappa)^{m_1}(\rho_1-\overline\kappa)^{m_1}=
(\rho_2-\kappa)^{m_1}(\rho_2-\overline\kappa)^{m_1}=(\rho_3-\kappa)^{m_1}(\rho_3-\overline\kappa)^{m_1},\\
&(\rho_1-\kappa)^{m_2}(\rho_1-\overline\kappa)^{m_2}(\rho_1^2+1)^{m_3-m_2}=
(\rho_2-\kappa)^{m_2}(\rho_2-\overline\kappa)^{m_2}(\rho_2^2+1)^{m_3-m_2}\\
&\quad=
(\rho_3-\kappa)^{m_2}(\rho_3-\overline\kappa)^{m_2}(\rho_3^2+1)^{m_3-m_2}\end{cases}$$
which leads to an algebraic relation on $\kappa$ and
$\overline\kappa$, contradicting their generality. \end{proof}

\begin{lemma}\label{lnn9}
In the notation of section \ref{secn2}, let $\Del_k$ be a triangle
in $S_C$, dual to a trivalent vertex of the graph $G$. Then $C_k$
is the union of two imaginary conjugate nodal rational curves,
which cross each other transversally at their non-singular points.
Furthermore, each of these rational curves crosses
$\Tor(\partial\Del_k)$ at precisely three imaginary points, where
it is non-singular. The number of real solitary nodes of $C_k$ is
equal to $|\Int(\Del_k)\cap\Z^2|+1$ $\mod 2$. Given the
intersection points of $C_k$ with $\Tor(\sig_1)$, $\Tor(\sig_2)$,
$\sig_1,\sig_2$ being two edges of $\Del_k$, the number of curves
like $C_k$ is equal to $2|\Del_k|(|\sig_1|\cdot|\sig_2|)^{-1}$.
\end{lemma}

\begin{proof} From the proof of Proposition \ref{pn1} we know that $C_k$
must be the union of two conjugate rational curves, each of them
crossing $\Tor(\partial\Del_k)$ at three points, where they are
unibranch. Again along the proof of Proposition \ref{pn1}, the two
pairs of conjugate intersection points of $C_k$ with
$\Tor(\sig_1)\cup\Tor(\sig_2)$ are generic, and thus, the
components of $C_k$ are nodal and non-singular along
$\Tor(\partial\Del_k)$, they intersect transversally and only at
their non-singular points. In particular, $C_k$ has no
non-solitary real nodes, and hence the number of its solitary real
nodes has parity of $|\Int(\Del_k)\cap\Z^2|+1$. For a component of
$C_k$, we have two choices of a pair of intersection points with
$\Tor(\sig_1)\cup\Tor(\sig_2)$ up to conjugation. Furthermore,
fixing these intersection points, we have by \cite{ShP}, Lemma
3.5,
$(|\Del_k|/4)/(|\sig_1|/2\cdot|sig_2|/2)=|\Del_k|/(|\sig_1|\cdot|\sig_2|)$
choices for a given component, and the proclaimed number of
choices for $C_k$ follows. \end{proof}

\subsection{Refinement of the tropical limit}\label{secnn4} Using the structure
of the curves $C_1,...,C_N$ which appear in the tropical limit of
the curve $C\subset\Tor_\K(\Del)$, we refine the tropical limit
along the procedures described in \cite{ShP}, sections 3.5 and
3.6. Namely, the local branches of the curves $C_1,...,C_N$, which
are tangent to the divisors
$\Tor(\partial\Del_1),...,\Tor(\partial\Del_N)$, are naturally
combined into disjoint pairs $(B,B')$ such that either
\begin{itemize}\item $B$ and $B'$ are centered at the same
point $z\in\Tor(\sig)$, where $\sig$ is a common edge of two
triangles, or \item $B$ is centered at a point $z\in\Tor(\sig)$;
$B'$ is centered at a point $z'\in\Tor(\sig')$ with $\sig,\sig'$
being parallel edges of two triangles, dual to the same extended
edge of $A_C$, and such that $z$ and $z'$ are joined by a sequence
of components of $C^{(0)}$ from the set ${\mathcal C}_b$.
\end{itemize}
Along the procedures of \cite{ShP}, sections 3.5 and 3.6, for any
pair $(B,B')$, the polynomial $f(x,y)$, defining the curve $C$,
determines a complex polynomial $f_{B,B'}(x,y)$ such that (i) it
has Newton triangle $\Del_{B,B'}=\conv\{(0,2),(0,0),(m,1)\}$,
where $m$ is the intersection number of $B$ with $\Tor(\sig)$, and
the coefficient of $x^{m-1}y$ is zero, (ii) the polynomial
$f_{B,B'}$ defies a rational nodal curve $C_{B,B'}$ in
$\Tor(\Del_{B,B'})$. Furthermore, the polynomials $f_1,...,f_N$,
defining the curves $C_1,...,C_N$ and given by (\ref{enn7}),
determine the coefficients of $f_{B,B'}$ at the vertices of
$\Del_{B,B'}$.

Among the results of \cite{ShP}, Lemma 3.9 and Proposition 6.1, we
shall use the following one:

\begin{lemma}\label{lnn10} In the above notation, assume that $m\ge 2$
and $a,b,c\in\C^*$. Then there are precisely $m$ polynomials
\begin{equation}\varphi(x,y)=ay^2+y\left(bx^m+\sum_{i=0}^{m-2}b_ix^i\right)+c\ ,\label{enn35}\end{equation} which
define nodal rational curves in $\Tor(\Del_{B,B'})$. If $a,b,c$
are real and $m$ is even, then there are precisely two real
polynomials (\ref{enn35}), which define rational nodal curves in
$\Tor(\Del_{B,B'})$; furthermore, one of these curves has $m-1$
real solitary nodes, and the other has no such nodes. If $a,b,c$
are real and $m$ is odd, then there is only one real polynomial
(\ref{enn35}), which defines a rational nodal curve, and this
curve has $m-1$ real solitary nodes.
\end{lemma}

\subsection{The
Welschinger number of a real rational tropical curve}\label{secn3}
Let $A$ be a real rational tropical curve of type $(r',r'',s'')$
with Newton polygon $\Del$, passing through
$\overline\bx'=\{\bx'_1,...,\bx'_{r'}\}$ and
$\overline\bx''=\{\bx''_1,...,\bx''_{r''}\}$ in such a way that
$\bx''_1,...,\bx''_{s''}$ are trivalent vertices of $A$, and the
points $\bx'_1,...,\bx'_{r'}$ and $\bx''_{s''+1},...,\bx''_{r''}$
lie in the interior parts of some edges of $A$. We intend to
define the Welschinger number $w(A,\overline\bx',\overline\bx'')$,
which later will be equated with the contribution to the
Welschinger number of all the real rational curves in the linear
system $|D|$ on the surface $\Sig=\Tor_\K(\Del)$, which pass
through the points $\bp'_i$, $i=1,...,r'$, and
$\bp''_{i1},\bp''_{i2}$, $i=1,...,r''$, and project onto the
tropical curve $A$.

Among the possible graphs $G$, inscribed in the given tropical
curve $A$ along the rules presented at the end of section
\ref{secn2}, there is the maximal graph $G_{\max}$. It is uniquely
determined by the disposition of the sets $\overline\bx'$ and
$\overline\bx''$ on $A$, and is obtained by letting
$V_0=V^2(G_0)$, $V_1=V^2(G_1)$, and so on, in the notation of
section \ref{secn2}.

If there is an even weight edge of $A$ lying outside $G_{\max}$,
we put $w(A,\overline\bx',\overline\bx'')=0$.

Assume that the edges of $A$ outside $G_{\max}$ have odd weights.
Let $S$ be the dual subdivision of $\Del$. Then the triangles of
$S$ dual to the vertices $\bx''_1,...,\bx''_{s''}$ of $A$ have all
edges of odd length, the triangles dual to the univalent vertices
of $G_{\max}$ have two odd length edges and one edge of even
length, the triangles dual to the trivalent vertices of $G_{\max}$
are {\it even}, that is, all edges are of even length. Then we put
\begin{equation}w(A,\overline\bx',\overline\bx'')=(-1)^{a+b}2^{-c}\prod|\Del'|\
,\label{enn100}\end{equation} where $a$ is the number of integral
points, which lie strictly inside the triangles of $S$, $b$ is the
number of the even triangles, $c$ is the number of of the triangles
of even lattice area, and $\Del'$ runs over all the triangles of
even lattice area, or dual to $\bx''_1,...,\bx''_{s''}$.

\section{The correspondence theorem}\label{secn14}

\subsection{Calculation of the Welschinger number via tropical
curves}\label{secn6} Let $\Sig$ be a Del Pezzo toric surface,
associated with a lattice polygon $\Del$ as defined in section
\ref{sec1}, and $D\subset\Sig$ be the corresponding real ample
divisor. Let the non-negative integers $r',r''$ satisfy
(\ref{enn13}).

\begin{theorem}\label{tnn1} In the above notation and assumptions,
let a collection of distinct points
$\bx'_1,...,\bx'_{r'},\bx''_1,...,\bx''_{r''}$ in $Q^2$ be
$\Del$-generic, and the distinct points $\bp'_i\in(\K_\R^*)^2$,
$i=1,...,r'$, and
$\bp''_{j1},\bp''_{j2}\in(\K^*)^2\backslash(\K_\R^*)^2$,
$j=1,...,r''$, be generic among those, satisfying
$$\val(\bp'_i)=\bx'_i,\ i=1,...,r',\quad\val(\bp''_{j1})=\val(\bp''_{j2})=\bx''_j,\ \conj(\bp''_{j1})=\bp''_{j2},
\quad j=1,...,r''\ .$$ Then one has
$$W_{r''}(\Sig,D)=\sum_{s''=0}^{r''}\sum_{A^{(s'')}}w(A^{(s'')},\overline\bx',
\overline\bx'')\ ,$$ where
$\overline\bx'=\{\bx'_1,...,\bx'_{r'}\}$,
$\overline\bx''=\{\bx''_1,...,\bx''_{r''}\}$, and $A^{(s'')}$
ranges over all real rational tropical curves of type
$(r',r'',s'')$ with Newton polygon $\Del$, which pass through
$\overline\bx'\cup\overline\bx''$ in such a way that precisely
$s''$ points of $\overline\bx''$ are trivalent vertices of
$A^{(s'')}$, and the other points of $\overline\bx''$ are interior
points of edges of $A^{(s'')}$, having even weight.
\end{theorem}

\begin{proof} {\it Step 1}. A $\Del$-general position of the points
$\overline\bx'\cup\overline\bx''$ and the general choice of the
respective points $\bp'_i$, $i=1,...,r'$, and
$\bp''_{j1},\bp''_{j2}$, $j=1,...,r''$, ensure all the generality
assumptions used in the proof of Proposition \ref{pn1}. Hence, in
particular, the real rational curves in the linear system $|D|$ on
the surface $\Tor_\K(\Del)$, passing through $\bp'_i$,
$i=1,...,r'$, and $\bp''_{i1},\bp''_{i2}$, $i=1,...,r''$, are
projected by $\val:(\K^*)^2\to\R^2$ to the real rational tropical
curves as described in the assertion of Theorem \ref{tnn1}.

\medskip
{\it Step 2}. Let $A$ be a real rational tropical curve of type
$(r',r'',s'')$, satisfying the conditions of Theorem \ref{tnn1},
$S$ the dual subdivision of $\Del$ into polygons
$\Del_1,...,\Del_N$, and $G$ a graph constructed as described at
the end of section \ref{secn2}. Without loss of generality assume
that $\bx''_1,...\bx''_{s''}$ are vertices of $A$. A collection of
real polynomials $f_1(x,y),...,f_N(x,y)$ with Newton polygons
$\Del_1,...,\Del_N$, respectively, is called an {\it admissible
tropicalization} if
\begin{itemize}\item for any common edge $\sig=\Del_k\cap\Del_l$, the truncations
$f_k^\sig$ and $f_l^\sig$ coincide; \item if an edge
$\sig\subset\Del_k$ is dual to an edge of $A$, passing through
$\bx'_i$, $1\le i\le r'$ (or through $\bx''_i$, $i>s''$), then
$f_k^\sig(\xi'_i,\eta'_i)=0$ (resp.,
$f_k^\sig(\xi''_i,\eta''_i)=f_k^\sig(\overline{\xi''_i},\overline{\eta''_i})=0$);
\item if a triangle $\Del_k$ is dual to $\bx''_i$,
$1\le i\le s''$, then
$f_k(\xi''_i,\eta''_i)=f_k(\overline{\xi''_i},\overline{\eta''_i})=0$;
\item if $\Del_k$ is a parallelogram, then $f_k$ is split into the
product of a monomial and few binomials; \item if a triangle
$\Del_k$ is dual to a vertex of $A$ outside the graph $G$, then
$f_k$ defines a real rational curve in $\Tor(\Del_k)$, which
crosses $\Tor(\partial\Del_k)$ at precisely three points and is
non-singular there; \item if a triangle $\Del_k$ is dual to a
univalent (or bivalent) vertex of $G$, then $f_k$ defines a real
rational curve in $\Tor(\Del_k)$, which crosses
$\Tor(\partial\Del_k)$ at two real and a pair of imaginary
conjugate points (resp., at one real and two pairs of imaginary
conjugate points) and is non-singular there; \item if a triangle
$\Del_k$ is dual to a trivalent vertex of $G$, then the curve
$f_k=0$ in $\Tor(\Del_k)$ is split into a pair of imaginary
conjugate rational curves, each of them crossing
$\Tor(\partial\Del_k)$ at precisely three imaginary points being
non-singular there. \end{itemize}

Next we define an {\it admissible refined tropicalization}
associated with $A$ and $G$. First, we consider the set of the
local branches of curves $C_k:=\{f_k=0\}\subset\Tor(\Del_k)$,
$k=1,...,N$, which are tangent to $\bigcup_k\Tor(\partial\Del_k)$,
and, following the instructions of section \ref{secnn4},
distribute these branches into disjoint pairs. With each pair
$(B,B')$ we associate a triangle
$\Del_{B,B'}=\conv\{(0,0),(0,2),(m,1)\}$, where $m$ is the
intersection number of $B$ (or $B'$) with
$\bigcup_k\Tor(\partial\Del_k)$. If $B$ and $B'$ are branches of
curves $C_k$ and $C_l$, respectively, centered at the same point
on $\Tor(\sig)$, $\sig=\Del_k\cap\Del_l$, we perform the
refinement procedure described in \cite{ShP}, section 3.5. The
corresponding transformations of the polygons $\Del_k,\Del_l$ and
the polynomials $f_k,f_l$, namely, the coordinate change induced
by $M_\sig\in\Aff(\Z^2)$ and the shift $x\mapsto x+\xi$, are
illustrated in Figure \ref{fign2} (see \cite{ShP}, section 3.5,
for all details). Notice that the (transformed) polynomials
$f_k,f_l$ determine the coefficients of $y^2$, $yx^m$, and $1$.

\begin{figure}
\setlength{\unitlength}{1cm}
\begin{picture}(13,6)(0,0)
\thinlines \put(0,0){\vector(1,0){3}}\put(0,0){\vector(0,1){6}}
\thicklines\put(1,2){\line(1,-1){1}}\put(2,1){\line(1,2){1}}
\put(1,2){\line(1,3){1}}\put(1,2){\line(2,1){2}}\put(2,5){\line(1,-2){1}}
\thinlines\put(3.5,2.9){\line(1,0){1}}\put(3.5,3.1){\line(1,0){1}}
\put(4.6,3.0){\line(-1,1){0.4}}\put(4.6,3.0){\line(-1,-1){0.4}}
\put(5,2){\vector(1,0){3}}\put(5,0){\vector(0,1){6}}
\thicklines\put(5,3){\line(1,0){2}}\put(5,3){\line(1,2){1}}
\put(5,3){\line(1,-2){1}}\put(6,5){\line(1,-2){1}}\put(6,1){\line(1,2){1}}
\thinlines\put(8.5,2.9){\line(1,0){1.5}}\put(8.5,3.1){\line(1,0){1.5}}
\put(10.1,3.0){\line(-1,1){0.4}}\put(10.1,3.0){\line(-1,-1){0.4}}\put(4.7,3){$1$}
\put(11,2){\vector(1,0){3}}\put(11,0){\vector(0,1){6}}\dashline{0.2}(11,3)(13,3)
\dashline{0.2}(13,3)(13,2)
\thicklines\put(11,4){\line(2,-1){2}}\put(11,5){\line(1,0){1}}
\put(12,5){\line(1,-2){1}}\put(11,2){\line(2,1){2}}\put(11,2){\line(0,-1){1}}
\put(11,1){\line(1,0){1}}\put(12,1){\line(1,2){1}}
\put(0.7,3.5){$\Del_k$}\put(1.2,1.0){$\Del_l$}\put(1.8,2.7){$\sig$}
\put(6.8,4){$\Del'_k$}\put(6.5,1.2){$\Del'_l$}\put(6,3.1){$\sig'$}
\put(12.8,4){$\Del''_k$}\put(11.7,0.5){$\Del''_l$}\put(12.9,1.6){$m$}
\put(10.6,3.9){$2$}\put(3.7,3.6){$M_\sig$} \put(8.7,3.3){\text{\rm
shift}}\put(10.7,2.8){$1$}
\end{picture}
\caption{Refinement of the tropicalization, I}\label{fign2}
\end{figure}
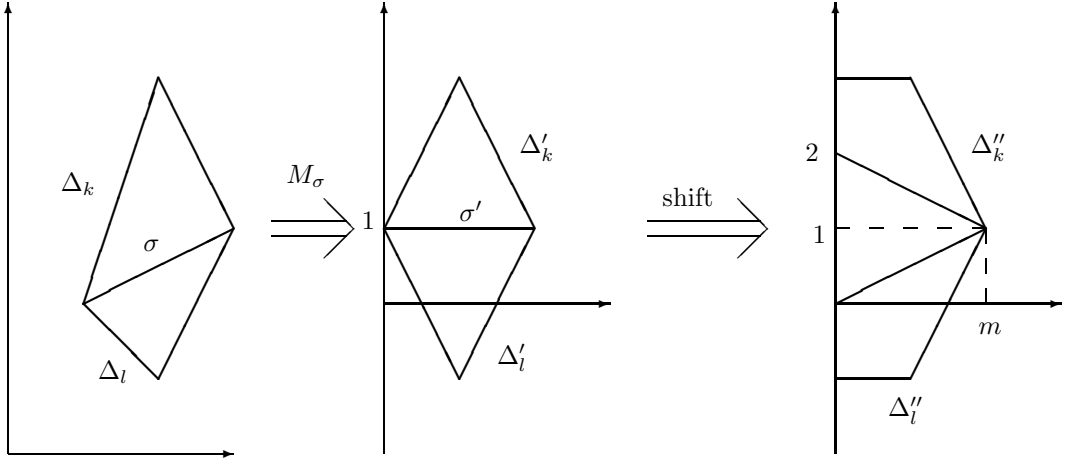

If branches $B,B'$, forming a pair, are centered at points
$z\in\Tor(\sig)$, $z'\in\Tor(\sig')$, where $\sig\ne\sig'$, then
we perform the refinement procedure, described in \cite{ShP},
section 3.6, and illustrated in Figure \ref{fign3}. Namely,
Figures \ref{fign3}(b,c) indicate the transforms induced by some
$M\in\Aff(\Z^2)$ and by a shift $x\mapsto x+\xi$. The convex
piece-wise linear function $\nu:\Del\to\R$, whose graph projects
onto the subdivision $S$, induces a convex piece-wise linear
function $\nu'$ on the polygons, shown in Figure \ref{fign3}(c)
except for the trapezoid $\theta$. The generality of the
configuration $\overline\bx'\cup\overline\bx''$ yields, in
particular, the generality of the function $\nu$ in the following
sense: no four parallel edges of the graph of $\nu$ lie on two
parallel planes. Thus, $\nu$ can uniquely be extended to $\theta$
as a convex piece-wise linear function, which defines a
subdivision of $\theta$ into parallelograms and a translate of
$\Del_{B,B'}$ (see, for example, Figure \ref{fign3}(d)).
Furthermore, those polynomials $f_1,...,f_N$, whose Newton
polygons contain $\sig,\sig'$ and all the parallel to them edges
of $S$ (see Figure \ref{fign3}(a)), after transformations
determine the coefficients at the integral points on the part of
the boundary of $\theta$ common with other polygons. This uniquely
determines the polynomials with Newton parallelograms inside
$\theta$, which are split into products of a monomial and
binomials. Finally, we determine the (non-zero) coefficients
corresponding to the vertices of the triangle inside $\theta$, and
which we respectively assign to the vertices of $\Del_{B,B'}$.

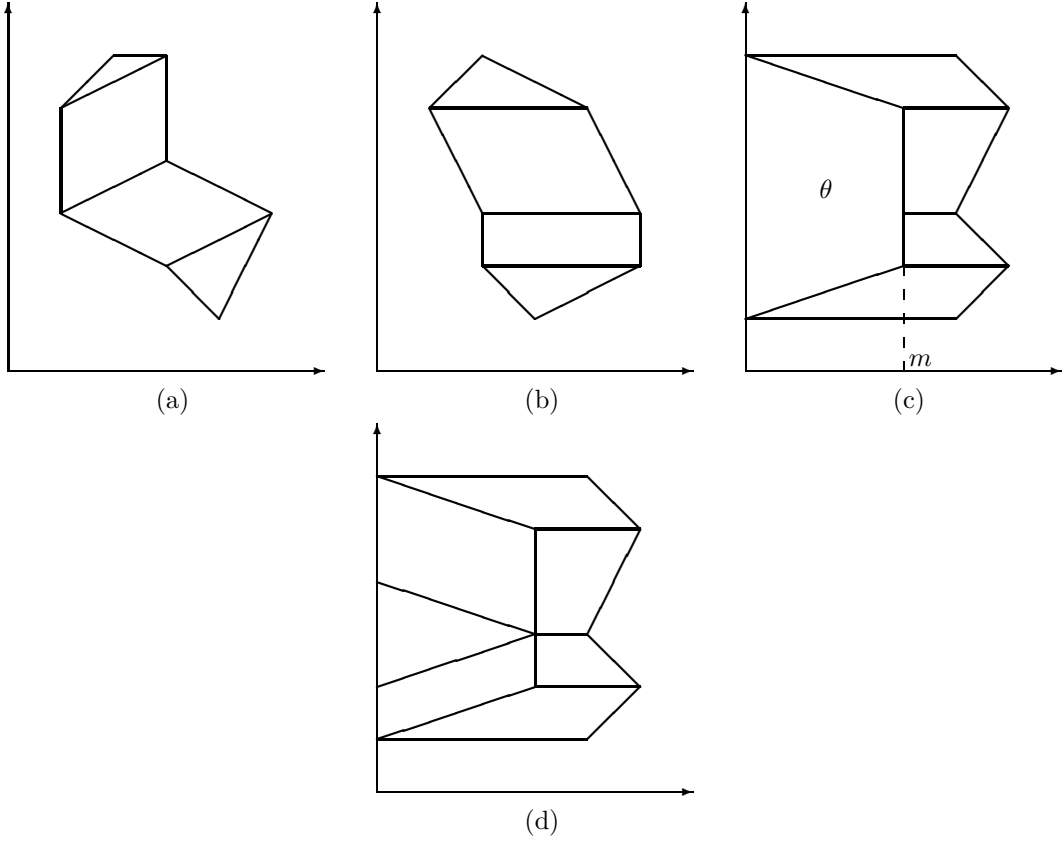
\begin{figure}
\setlength{\unitlength}{0.7cm}
\begin{picture}(14,16)(3,0)
\thinlines \put(0,9){\vector(1,0){6}}\put(0,9){\vector(0,1){7}}
\thicklines\put(1,12){\line(0,1){2}}\put(1,12){\line(2,1){2}}
\put(1,12){\line(2,-1){2}}\put(1,14){\line(1,1){1}}\put(1,14){\line(2,1){2}}
\put(2,15){\line(1,0){1}}\put(3,11){\line(1,-1){1}}
\put(3,11){\line(2,1){2}}\put(3,13){\line(0,1){2}}\put(3,13){\line(2,-1){2}}
\put(4,10){\line(1,2){1}}\put(2.8,8.3){\text{\rm (a)}} \thinlines
\put(7,9){\vector(1,0){6}}\put(7,9){\vector(0,1){7}}
\thicklines\put(8,14){\line(1,-2){1}}\put(8,14){\line(1,1){1}}
\put(9,11){\line(0,1){1}}\put(9,11){\line(1,0){3}}\put(8,14){\line(1,0){3}}
\put(9,11){\line(1,-1){1}}\put(9,12){\line(1,0){3}}
\put(9,15){\line(2,-1){2}}\put(10,10){\line(2,1){2}}\put(11,14){\line(1,-2){1}}
\put(12,11){\line(0,1){1}}\put(9.8,8.3){\text{\rm (b)}} \thinlines
\put(14,9){\vector(1,0){6}}\put(14,9){\vector(0,1){7}}
\thicklines\put(14,10){\line(3,1){3}}\put(14,10){\line(1,0){4}}
\put(18,10){\line(1,1){1}}\put(17,11){\line(1,0){2}}\put(17,11){\line(0,1){1}}
\put(17,12){\line(1,0){1}}\put(17,12){\line(0,1){2}}
\put(18,12){\line(1,-1){1}}\put(17,14){\line(1,0){2}}\put(18,12){\line(1,2){1}}
\put(14,15){\line(3,-1){3}}\put(14,15){\line(1,0){4}}\put(18,15){\line(1,-1){1}}\put(16.8,8.3){\text{\rm
(c)}}\thinlines\dashline{0.2}(17,9)(17,11)\put(17.1,9.1){$m$}\put(15.4,12.3){$\theta$}
\thinlines\put(7,1){\vector(1,0){6}}\put(7,1){\vector(0,1){7}}
\thicklines\put(7,2){\line(3,1){3}}\put(7,2){\line(1,0){4}}
\put(11,2){\line(1,1){1}}\put(10,3){\line(1,0){2}}\put(10,3){\line(0,1){1}}
\put(10,4){\line(1,0){1}}\put(10,4){\line(0,1){2}}
\put(11,4){\line(1,-1){1}}\put(10,6){\line(1,0){2}}\put(11,4){\line(1,2){1}}
\put(7,7){\line(3,-1){3}}\put(7,7){\line(1,0){4}}\put(11,7){\line(1,-1){1}}
\put(7,3){\line(3,1){3}}\put(7,5){\line(3,-1){3}}\put(9.8,0.3){\text{\rm
(d)}}
\end{picture}
\caption{Refinement of the tropicalization, II}\label{fign3}
\end{figure}

Denote by ${\mathcal B}$ the set of all pairs $(B,B')$ of local
branches as above. To obtain an admissible refined
tropicalization, we extend an admissible tropicalization
$f_1,...,f_N$ by adding (complex) polynomials $f_{B,B'}(x,y)$ with
Newton polygons $\Del_{B,B'}$, $(B,B')\in{\mathcal B}$, such that
\begin{itemize}\item any $f_{B,B'}$ is given by formula
(\ref{enn35}) with the coefficients $a,b,c\in\C^*$, prescribed as
explained above; \item the curve
$\{f_{B,B'}=0\}\subset\Tor(\Del_{B,B'})$ is rational nodal; \item
if $(B,B')$ is real then $f_{B,B'}$ is real; \item if $(B,B')$ is
imaginary, and $(\overline B,\overline B')$ is the conjugate pair,
then $f_{\overline B,\overline B'}=\overline f_{B,B'}$.
\end{itemize}

The curves $C_k=\{f_k=0\}\subset\Tor(\Del_k)$, $k=1,...,N$, and
$C_{B,B'}=\{f_{B,B'}=0\}\subset\Tor(\Del_{B,B'})$,
$(B,B')\in{\mathcal B}$, are called {\it admissible
tropicalization curves}.

\medskip
{\it Step 3}. Let us be given a tropical curve $A$ and an
admissible refined tropicalization $f_i$, $i=1,...,N$, $f_{B,B'}$,
$(B,B')\in{\mathcal B}$, as described in Step 2. We state

\begin{lemma}\label{lnn11}
(1) Under the above assumptions, there exists a polynomial
$f\in\K_\R[x,y]$ with Newton polygon $\Del$ such that
$\{f(x,y)=0\}$ is a rational nodal curve $C$ in $\Tor_\K(\Del)$,
and the refined tropical limit of the polynomial $f(x,y)$ consists
of the tropical curve $A$ and the polynomials $f_i$, $i=1,...,N$,
$f_{B,B'}$, $(B,B')\in{\mathcal B}$. All these curves $C$ have the
same number of real solitary nodes, which is equal to the total
number of such nodes over all the curves
$C_i=\{f_i=0\}\subset\Tor(\Del_i)$, $i=1,...,N$, and
$C_{B,B'}=\{f_{B,B'}=0\}\subset\Tor(\Del_{B,B'})$,
$(B,B')\in{\mathcal B}$.

(2) The number of the above curves $C$, passing through $\bp'_i$,
$i=1,...,r'$, and $\bp''_{i1},\bp''_{i2}$, $i=1,...,r''$, is
finite, and depends only on $A$ and on the polynomials $f_i$,
$i=1,...,N$. Furthermore, if all the edges of $A$ outside the
graph $G_{\max}$ have odd weight, then the number of such curves
$C$ is equal to $2^{s''-r''}\prod_{j>s''}|\sig_j|$, where $\sig_j$
is the edge of $S$, dual to the edge of $A$ passing through
$\bx''_j$, $j>s''$. \end{lemma}

The proof, based on \cite{ShP}, Theorem 5, is presented in section
\ref{secn5}.

\medskip
{\it Step 4}. By Proposition \ref{pn1}, the real rational curves
in the linear system $|D|$ on $\Tor_\K(\Del)$, passing through
$\bp'_i$, $i=1,...,r'$, and $\bp''_{i1},\bp''_{i2}$,
$i=1,...,r''$, project onto real rational tropical curves as
described in the assertion of Theorem \ref{tnn1}, and their
refined tropical limits determine (up to a real constant factor)
an admissible refined tropicalization, associated with $A$ and a
respective graph $G$.

We claim that, if there is an edge of $A$ of even weight, which is
not included in $G$, then the total contribution to
$W_{r''}(\Sig,D)$ of the real rational curves in $\Tor_\K(\Del)$,
having $A$ and $G$ in their tropical limits, is zero.

Indeed, by restriction (R3) obtained in Step 4 of the proof of
Proposition \ref{pn1}, semi-infinite edges of $A$ of even weight
belong to $G$. Let $f_k$, $k=1,...,N$, $f_{B,B'}$,
$(B,B')\in{\mathcal B}$, be an admissible tropicalization
associated with $A$ and $G$. An edge of $A$ of even weight outside
$G$ must have two endpoints, and hence there is a pair
$(B_0,B'_0)\in{\mathcal B}$ of real local branches, having an even
intersection multiplicity with $\bigcup_k\Tor(\partial\Del_k)$. By
Lemma \ref{lnn10}, the polynomial $f_{B_0,B'_0}$ can be replaced
by another polynomial $\hat f_{B_0,B'_0}$ keeping the
admissibility property. By Lemma \ref{lnn11}, these admissible
refined tropical limits produce an equal number of real rational
curves in $|D|$, but they have different parity of the number of
real solitary nodes (see again Lemma \ref{lnn10}).

Thus, the only pairs $(A,G)$, which contribute to
$W_{r''}(\Sig,D)$ are $(A,G_{\max})$ such that all edges of $A$
outside $G_{\max}$ have odd weight.

\medskip
{\it Step 5}. In the previous notation and assumptions, let $A$
have no edge of even weight outside $G_{\max}$. By Lemmas
\ref{lnn5}, \ref{lnn6}, \ref{lnn7}, \ref{lnn9}, and \ref{lnn10},
all the real rational curves in $|D|$ arising along Lemma
\ref{lnn11} from admissible refined tropicalizations associated
with $A$ and $G_{\max}$, have the same number of real solitary
nodes $\mod 2$, and the parity is given by the sign of
$w(A,\overline\bx',\overline\bx'')$ as defined in section
\ref{secn3}. Thus, to complete the proof of Theorem \ref{tnn1}, we
have to count the number of admissible refined tropicalizations
associated with $A$ and $G_{\max}$ and considered up to
multiplication by a non-zero real constant.

First, by Lemma \ref{lnn6}, we have $\prod|\Del'|$ choices for
admissible curves with Newton triangles, dual to
$\bx''_1,...,\bx''_{s''}$, where $\Del'$ runs over all these
triangles.

Second, the above choice and the coordinates of the points
$\bp'_i$, $i=1,...,r'$, and the points $\bp''_{i1}$, $i>s''$,
determine the intersections of the possible admissible
tropicalization curves $C_1,...,C_N$ with $\bigcup\Tor(\sig')$,
where $\sig'$ runs over all edges of the subdivision $S$, dual to
the extended edges of $A$, passing through $\bx'_i$, $i=1,...,r'$,
and $\bx''_i$, $i=1,...,r''$. Let $E_0$ be the union of all these
extended edges of $A$.

The graph $E_0$ is a union of trees, whose vertices differing from
$\bx''_1,...,\bx''_{r''}$ have valency at most two. Notice that
the position of the remaining edges of $A$ is prescribed by $E_0$
(provided that the combinatorial type of the pair
$(A,\overline\bx'\cup\overline\bx'')$ is fixed). This yields that
the graph $E_0$ has a bivalent vertex $v_0$, and there is an
extended edge $\eps_1$ in $A\backslash E_0$ starting at $v_0$.
Notice that $\eps_1$ cannot join two vertices of $E_0$ of valency
two (this can be verified precisely as a similar statement on the
graphs $G_i$, $i\ge 0$ in Step 3 of the proof of Proposition
\ref{pn1}). Put $E_1=E_0\cup\eps_1$.

For an admissible tropicalization curve with Newton triangle
$\Del'$, dual to $v_0$, we know the intersections with
$\Tor(\sig')\cup\Tor(\sig'')$ with $\sig',\sig''$, being any two
edges of the Newton triangle. Then, by Lemmas \ref{lnn5},
\ref{lnn7}, and \ref{lnn9}, the choice of such an admissible
tropicalization curve can be made (i) in a unique way, if all the
edges of $\Del'$ have odd length, (ii) in $|\Del'|/|\sig'|$ ways,
if $\sig'$ has an even length and $\sig''$ has an odd length,
(iii) in $2|\Del'|/(|\sig'|\cdot|\sig''|)$ ways, if both $\sig'$
and $\sig''$ have even length.

Next, by a similar reason, there is a bivalent vertex $v_1$ of
$E_1$, and the extended edge $\eps_2$ of $A\backslash E_1$,
starting at $v_1$, which does not end up at another bivalent
vertex of $E_1$. We then reconstruct a compatible admissible
tropicalization curve with Newton triangle dual to $v_1$.
Proceeding in the same manner, we complete the reconstruction of
all the admissible tropicalization curves $C_1,...,C_N$. It
follows immediately that the reconstruction of $C_1,...,C_N$ can
be done in
$2^e\prod|\Del'|\cdot(\prod\omega(\eps')\cdot\prod\omega(\eps''))^{-1}$
ways, where $e$ is the number of even triangles in $S$, $\Del'$
runs over all triangles in $S$, containing an even length edge,
$\eps'$ runs over all the extended edges of $A$ of even weight
with two endpoints, $\eps''$ runs over all the extended edges of
$A$, passing through $\bx''_j$, $j>s''$, and at last, $\omega(*)$
denotes the weight.

By Lemma \ref{lnn10}, for any collection of admissible
tropicalization curves $C_1,...,C_N$, we can find
$\prod(\omega(\eps')/2)$ compatible collections of admissible
refining tropicalization curves $C_{B,B'}$, $(B,B')\in{\mathcal
B}$, where $\eps'$ runs over all the extended edges of $A$ of even
weight.

Combining these computations with the result of Lemma \ref{lnn11},
we complete the proof. \end{proof}

\subsection{Proof of Lemma \ref{lnn11}}\label{secn5}
{\it Step 1}. The subdivision $S$, the convex piece-wise linear
function $\nu:\Del\to\R$, and the collection of admissible
tropicalization curves $C_i$, $i=1,...,N$, $C_{B,B'}$,
$(B.B')\in{\mathcal B}$, satisfy the conditions of the
patchworking Theorem 5 from \cite{ShP}, section 5.3. Indeed, the
assumptions of \cite{ShP}, section 5.1, are easily verified. An
orientation $\Gamma$ of the tropical curve $A$, required in the
assertion of Theorem 5 \cite{ShP}, can be chosen so that all edges
of $A$ form angles $\alp\in(-\pi/2,\pi/2]$ with the horizontal
coordinate axis. Then all the triads
$(\Del_k,\Del_k^-(\Gamma),C_k)$, $k=1,...,N$, and all the
deformation patterns, given by the curves $C_{B,B'}$,
$(B,B')\in{\mathcal B}$, are transversal (in the sense of
\cite{ShP}, section 5.2) by \cite{ShP}, Lemmas 5.5 and 5.6. Hence
\cite{ShP}, Theorem 5, provides the existence of a real nodal
rational curve $C$ in the linear system $|D|$ on $\Tor_\K(\Del)$,
whose refined tropical limit consists of $A$ and $C_i$,
$i=1,...,N$, $C_{B,B'}$, $(B.B')\in{\mathcal B}$. Moreover, the
number of real solitary nodes of any of such curves $C$ is equal
to the total number of real solitary nodes of the refined
tropicalization $C_i$, $i=1,...,N$, and $C_{B,B'}$,
$(B,B')\in{\mathcal B}$.

To further impose the conditions $\bp'_i,\bp''_{j1},\bp''_{j2}\in
C$, $i=1,...,r'$, $j=1,...,r''$, we need a description of the
complete family of curves $C\subset\Tor_K(\Del)$ as above. Such a
description is not given in the assertion of Theorem 5 \cite{ShP},
but we shall extract it from the proof of that theorem, for which
we refer to \cite{ShPo}. The proof goes as follows. The curves $C$
are represented by polynomials
$$f(x,y)=\sum_{(i,j)\in\Del}(a_{ij}+c_{ij})t^{\nu(i,j)}x^iy^j\in\K_\R[x,y]\
,$$ where $$f_k(x,y)=\sum_{(i,j)\in\Del_k}a_{ij}x^iy^j,\
k=1,...,N,\quad c_{ij}=c_{ij}(t)\in\K_\R,\ \val(c_{ij})<0,\
(i,j)\in\Del\ .$$ The rationality requirement is expressed as a
system of equations on the unknowns $c_{ij}(t)$, whose
linearization at $t=0$ appears to be independent, and thus, the
implicit function theorem applies, resulting in the existence of a
family of the aforementioned curves $C\subset\Tor_\K(\Del)$,
smoothly depending on some parameters. Thus, we shall assume that
$c_{v_0}=0$, where $v_0$ is some vertex of $\Del$, since the
polynomials in $\K[x,y]$ are considered up to a constant factor,
and then we have only to indicate, which of $c_{ij}$,
$(i,j)\in\Del\backslash\{v_0\}$, can be chosen as independent
parameters.

\medskip {\it Step2}. We start by reducing the set of
variables $c_{ij}$, $(i,j)\in\Del\backslash\{v_0\}$ up to the
following set.

Denote by $V'(S)$ the set of (integral) middle points $v(\sig)$ of
all the edges $\sig$ of $S$, dual to the extended edges of $A$
belonging to the graph $G_{\max}$. We shall show that independent
variables can be chosen among
\begin{equation}c_{ij},\quad(i,j)\in V(S)\cup
V'(S)\backslash\{v_0\}\ ,\label{enn43}\end{equation} where $V(S)$
stands for the set of vertices of $S$.

Along the proof of Theorem 5 as presented in \cite{ShPo}, the
rationality condition for $C$ reduces to equations for $c_{ij}$,
$(i,j)\in\Del\backslash\{v_0\}$, which, up to terms containing $t$
to a positive power, coincide with the equations of the
transversality of all the triples $(\Del_k,\Del_k^-(\Gam),C_k)$,
$k=1,...,N$, which are written for the coefficients of the
polynomials $$\widetilde
f_k=\sum_{(i,j)\in\Del_k}(a_{ij}+c_{ij})x^iy^j,\quad
k=1,...,N,\quad c_{v_0}=0\ .$$

The transversality condition for a parallelogram $\Del_k$ (recall
that $\Del_k^-(\Gam)$ consists of two non-parallel edges) implies
by \cite{ShP}, Lemma 5.1, $$c_{ij}=L_{ij}\bigg(\{c_{\alp\bet}\ :\
(\alp,\bet)\in\Del_k^-(\Gam)\backslash\{v_0\}\}\bigg)$$
\begin{equation} +\ \text{h.o.t.}+O(t),\quad
(i,j)\in\Del_k\backslash\Del_k^-(\Gam)\
,\label{enn39}\end{equation} with some real linear functions
$L_{ij}$, $(i,j)\in\Del_k\backslash\Del_k^-(\Gam)$, as well as the
following relation: for any pair $\sig,\sig'$ of parallel edges of
$\Del_k$, it holds that
\begin{equation}\frac{a_{ij}+c_{ij}}{a_{i'j'}+c_{i'j'}}=\frac{a_{i_1j_1}+c_{i_1j_1}}{a_{i'_1j'_1}+
c_{i'_1j'_1}}+O(t),\quad (i,j),(i_1j_1)\in\sig,\
(i',j'),(i'_1,j'_1)\in\sig'\ ,\label{enn41}\end{equation} where
$(i,j)$ and $(i_1,j_1)$ run over the integral points of $\sig$,
and $(i',j'),(i'_1,j'_1)\in\sig'$ are obtained from
$(i,j),(i_1,j_1)$, respectively, by a shift taking $\sig$ to
$\sig'$.

By \cite{ShP}, Definition 5.2, the transversality of a triple
$(\Del_k,\Del_k^-(\Gam),C_k)$ with a triangle $\Del_k$ means that
the condition to pass through the nodes of the curve $C_k$, and
the relations
\begin{equation}\begin{cases}(C\cdot\Tor(\sig))_z\ge(C_k\cdot\Tor(\sig))_z,\quad z\in
C_k\cap\Tor(\sig),\quad \sig\subset\Del_k^-(\Gam),&\\
(C\cdot\Tor(\sig))_z\ge(C_k\cdot\Tor(\sig))_z-1,\quad z\in
C_k\cap\Tor(\sig),\quad
\sig\not\subset\Del_k^-(\Gam),&\end{cases}\label{enn40}\end{equation}
are independent for curves $C$ in the linear system $|C_k|$ on
$\Tor(\Del_k)$. From the description of the curves $C_k$ and the
Riemann-Roch theorem (in the form of Lemma 5.5 \cite{ShP}), it can
easily be verified that conditions (\ref{enn40}) result in
\begin{eqnarray}&c_{ij}=L_{ij}\big(\{c_{\alp\bet}\
:\ (\alp,\bet)\in\Del_k^-(\Gam)\cup((V(S)\cup
V'(S))\cap\partial\Del_k)\backslash\{v_0\}\}\big)\nonumber\\
&\qquad +\ \text{h.o.t.}+O(t)\ , \label{enn42}\\ &
\quad(i,j)\in\Del_k\backslash\big(\Del_k^-(\Gam)\cup((V(S)\cup
V'(S))\cap\partial\Del_k)\big)\ .\nonumber\end{eqnarray}

Combining (\ref{enn39}) and (\ref{enn42}), we obtain that the
variables $c_{ij}$, $(i,j)\in\Del\backslash(V(S)\cup V'(S))$, can
be expressed via variables (\ref{enn43}).

\medskip {\it Step 3}. We finally select an independent set of
variables among (\ref{enn43}).

In each parallelogram $\Del_k$ there is a unique vertex, which
does not belong to $\Del_k^-(\Gam)$. We remove all such vertices
from $V(S)$ and denote the complement set by $V_\Gam(S)$. Without
loss of generality, assume that $v_0\in V_\Gam(S)$.

The extended edge of $A$, passing through a point $\bx''_j$,
$j>s''$, is dual to one or a few parallel edges of $S$; we pick up
one of these edges of $S$, denote it by $\sig_j$, and introduce
the set $V''(S)=\{v(\sig_j)\ :\ s''<j\le r''\}\subset V'(S)$.

We claim that $c_{ij}$, $(i,j)\in V_\Gam(S)\cup
V''(S)\backslash\{v_0\}$, is a complete set of independent
variables, parameterizing the family of real rational curves $C$
as in the assertion of Lemma \ref{lnn11}.

Indeed, first, we notice that $\Gam$ defines a partial order on
the set of parallelograms of $S$, which we complete up to a linear
order, and then, following the chosen order, we can restore all
the variables $c_{ij}$, $(i,j)\in V(S)\backslash\{v_0\}$, by means
of relations (\ref{enn41}), restricted to the vertices of the
parallelograms.

Second, we notice that if an extended edge $\eps$ of $A$ is dual
to several (parallel) edges of $S$, say
$\sig^{(1)},...,\sig^{(k)}$, then we can find $c_{ij}$ with
$(i,j)=v(\sig^{(i)})$, $1\le i<k$, using relations (\ref{enn41})
and the knowledge of $c_{\alp\bet}$,
$(\alp,\bet)\in\{v(\sig^{(k)})\}\cup V(S)\backslash\{v_0\}$.

Third, having the variables $c_{\alp\bet}$, $(\alp,\bet)\in
V(S)\cup V''(S)\backslash\{v_0\}$, and using relations
(\ref{enn42}), we can find all the variables $c_{ij}$, $(i,j)\in
V'(S)\backslash V''(S)$. Namely, we follow the algorithm of
construction of the graph $G_{\max}$, starting with the graph
$G_0$ (see Steps 2 and 3 of the proof of Proposition \ref{pn1},
and the beginning of section \ref{secn3}). That is, we begin with
the knowledge of $c_{ij}$, $(i,j)=v(\sig)$ for all edges $\sig$ of
$S$, dual to the edges of the graph $G_0$. We successively append
extended edges of $A$ to $G_0$ so that each time we have a
triangle $\Del_k$ and the variables $c_{ij}$, where $(i,j)$ ranges
over the vertices of $\Del_k$ and over $v(\sig'),v(\sig'')$,
$\sig',\sig''$ being any two edges of $\Del_k$. Then relations
(\ref{enn42}) allow us to find $c_{ij}$ with $(i,j)=v(\sig''')$,
$\sig'''$ being the third edge of $\Del_k$.

Thus, the variables $c_{ij}$, $(i,j)\in V_\Gam(S)\cup
V''(S)\backslash\{v_0\}$, determine the variables $c_{ij}$,
$(i,j)\in V(S)\cup V'(S)\backslash\{v_0\}$, and hence, as shown
before, all the variables $c_{ij}$,
$(i,j)\in\Del_k\backslash\{v_0\}$. At last, we observe that
$\#(V_\Gam(S)\cup V''(S)\backslash\{v_0\})=\rk(A)+r''-s''=-K_\Sig
D-1$ is the dimension of the variety of rational curves in the
linear system $|D|$ on the surface $\Sig$, which confirms the
independence of the variables $c_{ij}$, $(i,j)\in V_\Gam(S)\cup
V''(S)\backslash\{v_0\}$.

A particular consequence of the result obtained in this step is
that the variables $c_{ij}$, $(i,j)\in V''(S)$ are independent of
the variables $c_{ij}$, $(i,j)\in V(S)\backslash\{v_0\}$, and that
equations (\ref{enn41}), restricted to the vertices of all the
parallelograms in $S$, impose $N_4$ independent conditions on the
latter set of variables\footnote{Recall that $N_4$ is the number
of parallelograms in $S$.}.

\medskip
{\it Step 4}. The requirement to pass through the points $\bp'_i$,
$i=1,...,r'$, and $\bp''_{i1},\bp''_{i2}$, $i=1,...,r''$, applied
to the family of the real rational curves $C$ with the given
refined tropical limit, can be formalized as follows.

Let a vertex $\bx''_i$, $1\le i\le s''$, of $A$ be dual to a
triangle $\Del_k$. Then the equations
$f(\bp''_{i1})=f(\bp''_{i2})=0$ coincide up to terms, containing
$t$ to a positive power, with equations
$f_k(\xi''_i,\eta''_i)=f_k(\overline{\xi''_i},\overline{\eta''_i})=0$
with $\xi''_i,\eta''_i$ introduced in (\ref{enn38}). By
Riemann-Roch (see Lemma \ref{lnn6}) the two latter equations
together with (\ref{enn42}), restricted to $\Del_k$, impose
$|\Del_k\cap\Z^2|-1$ independent conditions on the variables
$c_{\alp\bet}$, $(\alp,\bet)\in\Del_k\backslash\{v_0\}$, and thus,
one can conclude with
\begin{equation}\frac{a_{\alp_1\bet_1}+c_{\alp_1\bet_1}}{a_{\alp_3\bet_3}+c_{\alp_3\bet_3}}=O(t)\
,\quad
\frac{a_{\alp_2\bet_2}+c_{\alp_2\bet_2}}{a_{\alp_3\bet_3}+c_{\alp_3\bet_3}}=O(t)\
,\label{enn48}\end{equation} where
$(\alp_1,\bet_1),(\alp_2,\bet_2),(\alp_3,\bet_3)$ are the vertices
of $\Del_k$.

Let a point $\bx'_i$, $1\le i\le r'$, belong to an extended edge
$\eps$ of $A$, whose endpoints are dual to triangles
$\Del_k,\Del_l$ (the fact that $\eps$ has two endpoints is
established in Step 4 of the Proof of Proposition \ref{pn1}). Then
the condition $f(\bp'_i)=0$ can be expressed by an equation of
type (5.4.26) in \cite{ShP}, section 5.4, which, in our notation,
is interpreted as
\begin{eqnarray}&c_{\alp_1\bet_1}+\varphi_1(\xi'_i,\eta'_i)
c_{\alp_2\bet_2}\nonumber\\ &\qquad
=\varphi_2(\xi'_i,\eta'_i)\big(\varphi_3(\xi'_i,\eta'_i),a_{\alp_1\bet_1},
a_{\alp_2\bet_2},a'\big)^{1/m}t^p(1+\Phi)+\Psi\
,\label{enn46}\end{eqnarray} where $m$ is the length of the edge
$\sig$ of $\Del_k$, which is dual to $\eps$, the points
$(\alp_1,\bet_1)$ and $(\alp_2,\bet_2)$ are the vertices of
$\sig$, the non-zero real numbers $\xi'_i,\eta'_i$ are defined by
(\ref{enn37}), $\varphi_1,\varphi_2,\varphi_3$ are monomials of
the corresponding variables, $a'$ is a linear combination of
certain coefficients of $f_k$ (or $f_l$) such that $\widetilde a$
does not vanish\footnote{The coefficient $\widetilde a$ appears as
$a''_{01}$ in formula (5.4.26) of \cite{ShP}.}, $p$ is a positive
integer, and, at last, $\Phi$ and $\Psi$ are certain real analytic
functions of the variables $a_{ij}$, $c_{ij}$, $t$, whose terms
either contain $t$ to a positive power, or are at least quadratic
in $c_{ij}$. We observe that equation (\ref{enn46}) is split into
$m$ equations such that, for an odd $m$, only one of them is real,
and, for an even $m$, not one nor two of them are real.

Assume now that a point $\bx''_i$, $s''<i\le r''$, belongs to an
extended edge $\eps$ of $A$ with endpoints dual to triangles
$\Del_k$, $\Del_l$ (the case of a semi-infinite edge $\eps$ can be
treated in the same way). Again by formula (5.4.26) of \cite{ShP},
section 5.4, the conditions $f(\bp''_{i1})=f(\bp''_{i2})=0$ reduce
to a system of equations, which can be written as
\begin{equation}\begin{cases}&c_{\alp_1\bet_1}+\varphi_1(\xi''_i,\eta''_i)
c_{\alp_2\bet_2}\\ &\qquad
=\varphi_2(\xi''_i,\eta''_i)\bigg(\varphi_3(\xi''_i,\eta''_i),a_{\alp_1\bet_1},
a_{\alp_2\bet_2},\widetilde
a\bigg)^{1/m}t^p(1+\Phi)+\Psi,\\
&c_{\alp_1\bet_1}+\varphi_1(\overline{\xi''_i},\overline{\eta''_i})
c_{\alp_2\bet_2}\\
&\qquad=\varphi_2(\overline{\xi''_i},\overline{\eta''_i})
\bigg(\varphi_3(\overline{\xi''_i},\overline{\eta''_i}),a_{\alp_1\bet_1},
a_{\alp_2\bet_2},\overline
a'\bigg)^{1/m}t^p(1+\overline\Phi)+\overline\Psi,\end{cases}\label{enn47}\end{equation}
where $2m$ is the length of the edge $\sig$ of $\Del_k$, which is
dual to $\eps$, the point $(\alp_1,\bet_1)$ is a vertex of $\sig$,
different from $v_0$, and $(\alp_2,\bet_2)=v(\sig)$ is the middle
point of $\sig$, the non-zero numbers $\xi''_i,\eta''_i$ are
defined by (\ref{enn38}), $\varphi_1,\varphi_2,\varphi_3$ are
monomials of the corresponding variables with real coefficients,
$a'$ is a linear combination of certain coefficients of $f_k$ (or
$f_l$) such that $\widetilde a\ne 0$, $p$ is a positive integer,
and, at last, $\Phi$ and $\Psi$ are certain analytic functions of
the variables $a_{ij}$, $c_{ij}$, $t$, whose terms either contain
$t$ to a positive power, or are at least quadratic in $c_{ij}$. We
observe that system (\ref{enn47}) is split into $m$ systems
invariant with respect to the complex conjugation.

Take equations (\ref{enn41}), restricted to the vertices of all
the parallelograms in $S$, and equations (\ref{enn48}) for all
$i=1,...,s''$, and then append one real equation (\ref{enn46}) for
each $i=1,...,r'$, and one $\conj$-invariant system (\ref{enn47})
for each $i=s''+1,...,r''$. Observe that the linearizations of the
$|V(S)|-1+|V''(S)|$ chosen equations impose independent conditions
on the variables $c_{ij}$, $(i,j)\in V(S)\cup
V''(S)\backslash\{v_0\}$. Indeed, using these linearizations and
following the algorithm, which restores the tropical curve $A$
from the graph $E_0$ (see Step 5 in the proof of Theorem
\ref{tnn1}), we can determine all the aforementioned variables.

\section{Counting the Welschinger invariant via lattice paths in the
Newton polygon}\label{secn15}

\subsection{Reconstruction of tropical curve passing
through generic points on a straight line} Following the ideas of
\cite{M2,M3}, we intend to choose a specific configuration of the
points $\bx'_1,...,\bx'_{r'}$ and $\bx''_1,...,\bx''_{r''}$, for
which all the subdivisions $S$ of $\Del$, dual to the real
rational tropical curves mentioned in Theorem \ref{tnn1}, can be
obtained in a simple combinatorial algorithm. In this section we
describe a reconstruction of a real rational tropical curve.

We introduce an orthogonal system of coordinates $\lam=\alp x+\bet
y$, $\mu=\bet x-\alp y$ in $\R^2$ with generic $\alp,\bet\in\Q$.
The coordinate line $\Lam:=\{\mu=0\}$ is then not orthogonal to
any of the segments joining integral points in $\Del$.

We pick $r'+r''$ distinct points $\by_1,...,\by_{r'+r''}$ in
$\Lam\cap\Q^2$ one by one so that $0<\lam(\by_i)$ and
$\lam(\by_i)\ll\lam(\by_{i+1}$ for all $i\ge 1$.
\begin{equation}0<\lam(\by_i)\quad\text{and}\quad\lam(\by_i)\ll\lam(\by_{i+1})
\quad\text{for all}\quad i\ge 1\ ,\label{enn49}\end{equation} then
move them slightly in their neighborhoods in $\Q^2$, obtaining a
$\Del$-generic configuration $\overline\by$. At last, we declare
$r'$ points in $\overline\by$ as $\bx'_1,...,\bx'_{r'}$ and the
remaining $r''$ points as $\bx''_1,...,\bx''_{r''}$.

\begin{remark}
In fact, one could leave the set $\overline\by$ on the line
$\Lam$, and satisfy the generality conditions required in Theorem
\ref{tnn1}. In the construction of a real rational tropical curve
through $\overline\bx'\cup\overline\bx''$, presented below, we
suppose that the fixed points lie on the line $\Lam$, since a
small variation of the configuration does not affect the
construction.
\end{remark}

Let $A$ be a real rational tropical curve of type $(r',r'',s'')$
with Newton polygon $\Del$, which passes through
$\overline\bx'\cup\overline\bx''$ so that some $s''$ points among
$\bx''_1,...,\bx''_{r''}$ are its vertices. As established in Step
3 of the proof of Proposition \ref{pn1}, the other points in
$\overline\bx'\cup\overline\bx''$ are not vertices of $A$.

The curve $A$ can be restored in finitely many ways, and we shall
enumerate them.

\begin{figure}
\setlength{\unitlength}{1cm}
\begin{picture}(13,8)(0,0)
\thinlines\put(1,7){\line(1,0){11}}\thicklines\put(2,7){\vector(-1,1){1}}
\put(2,7){\vector(1,-1){1}}\put(4,7){\vector(0,1){1}}\put(4,7){\vector(0,-1){1}}
\put(7,7){\vector(-1,1){1}}\put(7,7){\vector(1,1){1}}\put(7,7){\vector(0,-1){1}}
\put(10,7){\vector(-1,-2){0.5}}\put(10,7){\vector(1,2){0.5}}\put(12.3,6.9){$\Lam$}
\put(1.8,6.5){$\bx'_1$}\put(7.2,6.5){$\bx''_1$}\put(7,5){$\text{\rm
(a)}$}\thinlines\put(1,2){\line(1,0){4.5}}\put(8,2){\line(1,0){4}}
\put(5.8,1.9){$\Lam$}\put(12.3,1.9){$\Lam$}\thicklines\put(2,2){\vector(1,1){2}}
\put(4,2){\vector(-1,1){2}}\put(9,2){\vector(1,1){1}}\put(11,2){\vector(-1,1){1}}
\put(10,3){\vector(0,1){1}}\put(5,1){$\text{\rm
(b)}$}\put(10,1){$\text{\rm (c)}$}
\end{picture}
\caption{Reconstruction of the tropical curve}\label{fign4}
\end{figure}
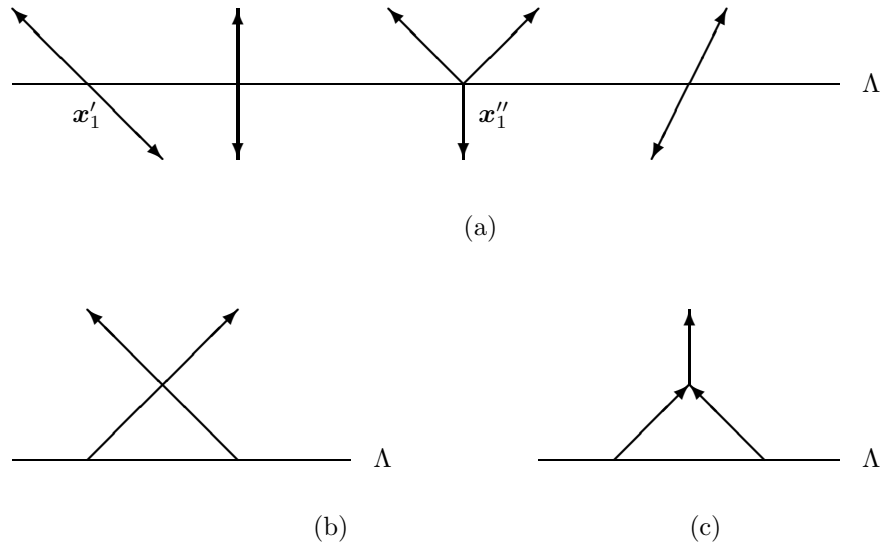

Consider the germs of edges of $A$ passing through
$\overline\bx'\cup\overline\bx''$. We have a finite choice of
these germs which should be orthogonal to some segments with
integral endpoints in $\Del$. Some $s''$ of the points
$\bx''_1,...\bx''_{r''}$ should be trivalent vertices of $A$. In
view of the general position of the line $\Lam$, the germs do not
lie on $\Lam$. Furthermore, we orient the germs or their halves to
start at the points $\bx'_1,...,\bx'_{r'},\bx''_1,...,\bx''_{r''}$
as shown, for example, on Figure \ref{fign4}(a) (cf. \cite{M3},
Lemma 4.17 and Figure 11).

Given a combinatorial type of $A$ and a combinatorial type of the
pair \mbox{$(A,\overline\bx'\cup\overline\bx'')$}, the
configuration $\overline\bx'\cup\overline\bx''$ determines $A$
uniquely. In particular, some extended edges of $A$ arising from
the given edge germs at $\overline\bx'\cup\overline\bx''$ must
intersect. For such an intersection point we have two
possibilities: (i) it is either a four-valent vertex of $A$ (see
Figure \ref{fign4}(b)), or (ii) it is a trivalent vertex of $A$
(see Figure \ref{fign4}(c)). In the latter case, another extended
edge of $A$ starts at the intersection point, which we orient as
shown in Figure \ref{fign4}(c). For this edge germ we again have a
finite choice, determined by the orthogonality to a segment with
integral endpoints in $\Del$ and by the equilibrium relation (see,
for instance, \cite{RST}, formula (10) in section 3), which, in
particular, says that the directing vector of the new edge germ is
a linear combination of the directing vectors of the two edges,
coming to the given vertex of $A$, with positive coefficients.

We proceed in the same manner with the new collection of edge
germs of $A$, and reconstruct the whole tropical curve $A$. In
fact, this procedure coincides with the reconstruction of $A$ in
the sequence of graphs $E_0, E_1,...$ as appears in Step 5 of the
proof of Theorem \ref{tnn1}.

We make a few observations: \begin{enumerate}\item[(O1)] The
directing vectors of the extended edges of $A$, which do not
intersect with $\overline\bx'\cup\overline\bx''$, and have an
endpoint in the half plane $\mu>0$ (or $\mu<0$), have a positive
(resp., negative) $\mu$-component, and hence $A\cap
\Lam=\overline\bx'\cup\overline\bx''$. This means that the
components of $\Lam\backslash(\overline\bx'\cup\overline\bx'')$
lie entirely in the components of $\R^2\backslash A$, and are dual
to some $r'+r''+1$ integral points in $\Del$ ordered by the linear
function $\lam$ in the same way as the components of
$\Lam\backslash(\overline\bx'\cup\overline\bx'')$. Furthermore,
the semi-infinite components of
$\Lam\backslash(\overline\bx'\cup\overline\bx'')$ correspond to
the (integral) points of $\Del$, where $\lam\big|_{\Del}$ takes
its maximal and minimal values.
\item[(O2)] If $\bx''_i$ is a trivalent vertex of $A$, and
$\omega_1,\omega_2,\omega_3$ are the integral points in $\Del$,
which correspond to the components of $\R^2\backslash A$ adjacent
to $\bx''_i$, so that the first two components intersect $\Lam$,
then $\lam(\omega_1)>\lam(\omega_3)>\lam(\omega_2)$.
\item[(O3)] Let ${\mathcal R}$ be the rays in the half-plane $\mu>0$,
starting at $\overline\bx'\cup\overline\bx''$ and at the trivalent
vertices of $A$, and generated by the oriented edge germs of $A$.
We introduce a partial order in ${\mathcal R}$ as follows. The
rays, starting at $\overline\bx'\cup\overline\bx''$, one orders by
successive intersections with a line close and parallel to $\Lam$
and oriented by $\lam$. If a ray $l_0\in{\mathcal R}$ starts at a
trivalent vertex, where rays $l',l''\in{\mathcal R}$ merge, and
$l'\succ l''$, then we assume $l_0\succ l$ for all $l\in{\mathcal
R}$ such that $l'\succ l$, and $l\succ l_0$ for all $l\in{\mathcal
R}$ such that $l\succ l'$. Then relation (\ref{enn49}) yields
that, if some rays $l_1,l_2\in{\mathcal R}$ cross a ray
$l_3\in{\mathcal R}$, then $l_1$ and $l_2$ are ordered, say,
$l_1\succ l_2$, and, respectively, $\mu(l_1\cap l_3)>\mu(l_2\cap
l_3)$. A similar claim holds for rays in the half-plane $\mu<0$.
\end{enumerate}

\subsection{Reconstruction of subdivisions dual to real rational tropical curves passing
through generic points on a straight line}\label{secn10}
\begin{figure}
\begin{center}
\epsfxsize 140mm \epsfbox{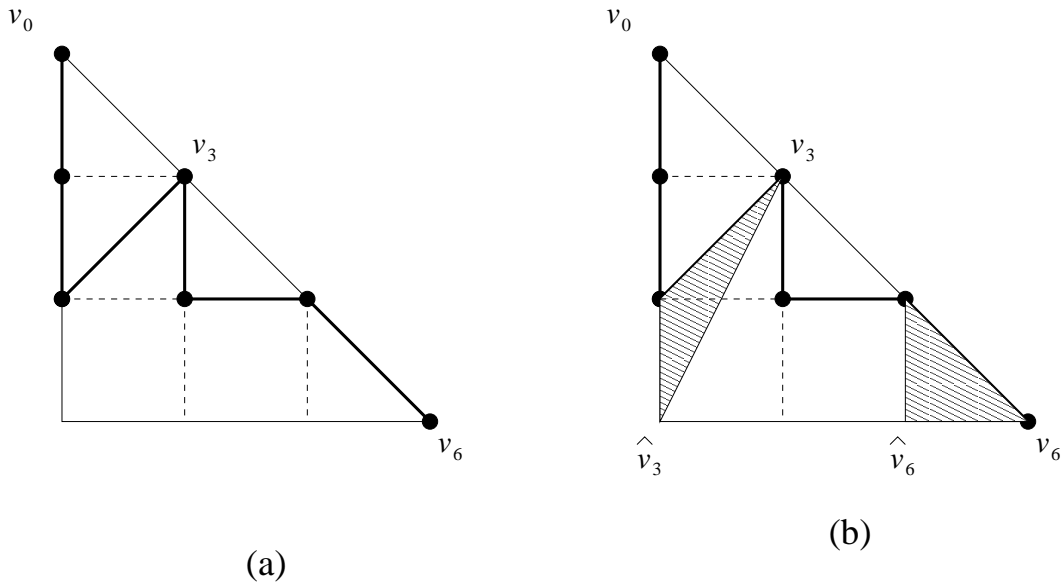}
\end{center}
\caption{Lattice path and subdivision $S_0$} \label{fign5}
\end{figure}

Observations O1, O2, O3 allow us to convert the preceding
construction of $A$ into the dual language of subdivisions of
$\Del$. The dual construction is performed as follows:

{\it Step 1}. The linear function $\lam$ defines a linear order of
the set $\Del\cap\Z^2$. We choose $r'+r''+1$ distinct successive
points $v_i\in\Del\cap\Z^2$, $i=0,...,r'+r''$, such that
$\lam(v_0)=\min\lam(\Del)$ and $\lam(v_{r'+r''})=\max\lam(\Del)$.
These points are dual to the components of $\R^2\backslash A$
containing the components of
$\Lam\backslash\overline\bx'\cup\overline\bx''$. We then connect
the points $v_0,...,v_{r'+r''}$ by segments
$\sig_i=[v_{i-1},v_i]$, $i=1,...,r'+r''$, obtaining a broken line,
called a {\it lattice path} (see Figure \ref{fign5}(a)).

{\it Step 2}. The segments $\sig_1,...,\sig_{r'+r''}$ naturally
correspond to the points $\by_1,...,\by_{r'+r''}\in \Lam$,
respectively. We take disjoint sets ${\mathcal V},{\mathcal
W}\subset\{1,2,...,r'+r''\}$ with $\#{\mathcal V}=s''$,
$\#{\mathcal W}=r'$, and declare \begin{itemize}\item the points
$\by_i$, $i\in{\mathcal W}$, to be (somehow ordered) points
$\bx'_1,...,\bx'_{r'}$, \item the points $\by_i$,
$i\not\in{\mathcal W}$, to be (somehow ordered) points
$\bx''_1,...,\bx''_{r''}$,
\item the points $\by_i$, $i\in{\mathcal V}$, to be trivalent vertices of $A$, respectively
dual to the lattice triangles in $\Del$.\end{itemize} That is, the
segments $\sig_j$, $j\not\in{\mathcal V}$, are dual to some edges
of $A$, whose germs pass through the corresponding points $\by_j$,
and are orthogonal to $\sig_j$, $j\not\in{\mathcal V}$,
respectively (see Figure \ref{fign4}(a)). In turn, if
$j\in{\mathcal V}$, then we pick a point $\hat
v_j\in\Del\cap\Z^2\backslash\sig_j$ such that
$\lam(v_{j-1})<\lam(\hat v_j)<\lam(v_j)$, and obtain a triangle
$T_j=\conv\{v_{j-1},\hat v_j,v_j)$ dual to $\by_j$ (see Figure
\ref{fign5}(b), where $T_j$ is designated by a shadow). The sides
of $T_j$ are dual to the edges of $A$, whose germs start at
$\by_j$ in the normal directions oriented outside of $T_j$.

In the absence of suitable points $\hat v_j$ for $j\in{\mathcal
V}$, or in the case of an odd length edge $\sig_i$ for
$i\not\in{\mathcal V}\cup{\mathcal W}$, we stop the construction
and say that the lattice path $\bigcup_i\sig_i$ and the sets
${\mathcal V}$, ${\mathcal W}$ are inconsistent.

{\it Step 3}. Put $S_0=\bigcup_i\sig_i\cup\bigcup_jT_j$. This is
part of a subdivision of $\Del$, and we extend it in the following
inductive process. Let $S_k$, $k\ge 0$, be a contractible union of
some lattice segments, triangles and parallelograms, such that
$S_k\supset S_0$, $S_k\ne\Del$, and, for each connected component
$\Del'$ of $\Del\backslash S_k$, the intersection
$\del=S_k\cap\partial\Del'$ is a lattice path in $\Del$ along
which the function $\lam$ is strongly monotone.

If $\conv(\del)\cap\Del'=\emptyset$ we stop the construction and
say that the sequence $S_0,...,S_k$ is inconsistent.

If $\conv(\del)\cap\Del'\ne\emptyset$, there a pair of successive
segments $\sig'$ and $\sig''$ in $\del$ such that
$\conv(\sig'\cup\sig'')\cap\Del'\ne\emptyset$, and
$\lam(\sig'\cap\sig'')$ is minimal among all such pairs. The pair
$\sig',\sig''$ corresponds to a pair of germs of edges of $A$,
generating rays which intersect each other, and, moreover, along
the partial order of rays as defined in (O3) above, this pair of
intersecting rays is minimal among all pairs of intersecting rays
which can be constructed in the current stage.

Then we define $S_{k+1}$ by adding to $S_k$ either the triangle
$\conv(\sig'\cup\sig'')$, or the parallelogram, built on
$\sig',\sig''$. For $A$ this means adding a new tri- or
four-valent vertex as shown in Figure \ref{fign4}(b,c).

If the polygon, added to $S_k$, is not contained in the closure of
$\Del'$, we stop the construction and call the sequence
$S_0,...,S_{k+1}$ inconsistent.

{\it Step 4}. Consider a sequence $S_0,...,S_k$ which ends up with
$S_k=\Del$. In parallel, the dual construction gives us a graph
$A\subset\R^2$. We assign corresponding weights to the edges of
$A$, and obtain a weighted rational graph, which by construction
satisfies the equilibrium condition at each vertex (see
\cite{RST}, formula (10) in section 3), and thus, by \cite{RST},
Theorem 3.6, is a tropical curve with Newton polygon $\Del$.

We call the obtained subdivision {\it consistent} if the tropical
curve $A$ is real rational of type $(r',r'',s'')$. For all the
consistent subdivisions of $\Del$, obtained in the above
algorithm, starting with the initial data
$v_0,...,v_{r'+r''},{\mathcal V},{\mathcal W}$, we sum up the
Welschinger numbers and denote the result by
$W_{\Del,r',r''}(v_0,...,v_{r'+r''}|{\mathcal V},{\mathcal W})$.
Moreover, we can get rid of the dependence on ${\mathcal V}$,
noticing that the non-zero Welschinger numbers come only from the
sets ${\mathcal V}$, defined as follows: $j\in{\mathcal V}$ if
$j\not\in{\mathcal W}$ and $|\sig_j|$ is odd.

Since all the real rational tropical curves passing through
$\overline\bx'\cup\overline\bx''$ are obtained in the above
construction, we conclude

\begin{theorem}\label{cnn1} Given a linear function
$\lam(x,y)=\alp x+\bet y$ with generic $\alp,\bet\in\Q$, and a set
${\mathcal W}\subset\{1,...,r'+r''\}$ with $\#{\mathcal W}=r'$,
one has
$$W_{r''}(\Sig,D)=\sum W_{\Del,r',r''}(v_0,...,v_{r'+r''}|{\mathcal W})\ ,$$
where the sum ranges over all sequences of points
$v_0,...,v_{r'+r''}\in\Del\cap\Z^2$ such that
$$\min\lam(\Del)=\lam(v_0)<...<\lam(v_{r'+r''})=\max\lam(\Del)\
.$$
\end{theorem}

\subsection{Example: real rational plane cubics}\label{secn12}
\begin{figure}
\begin{center}
\epsfysize 200mm \epsfbox{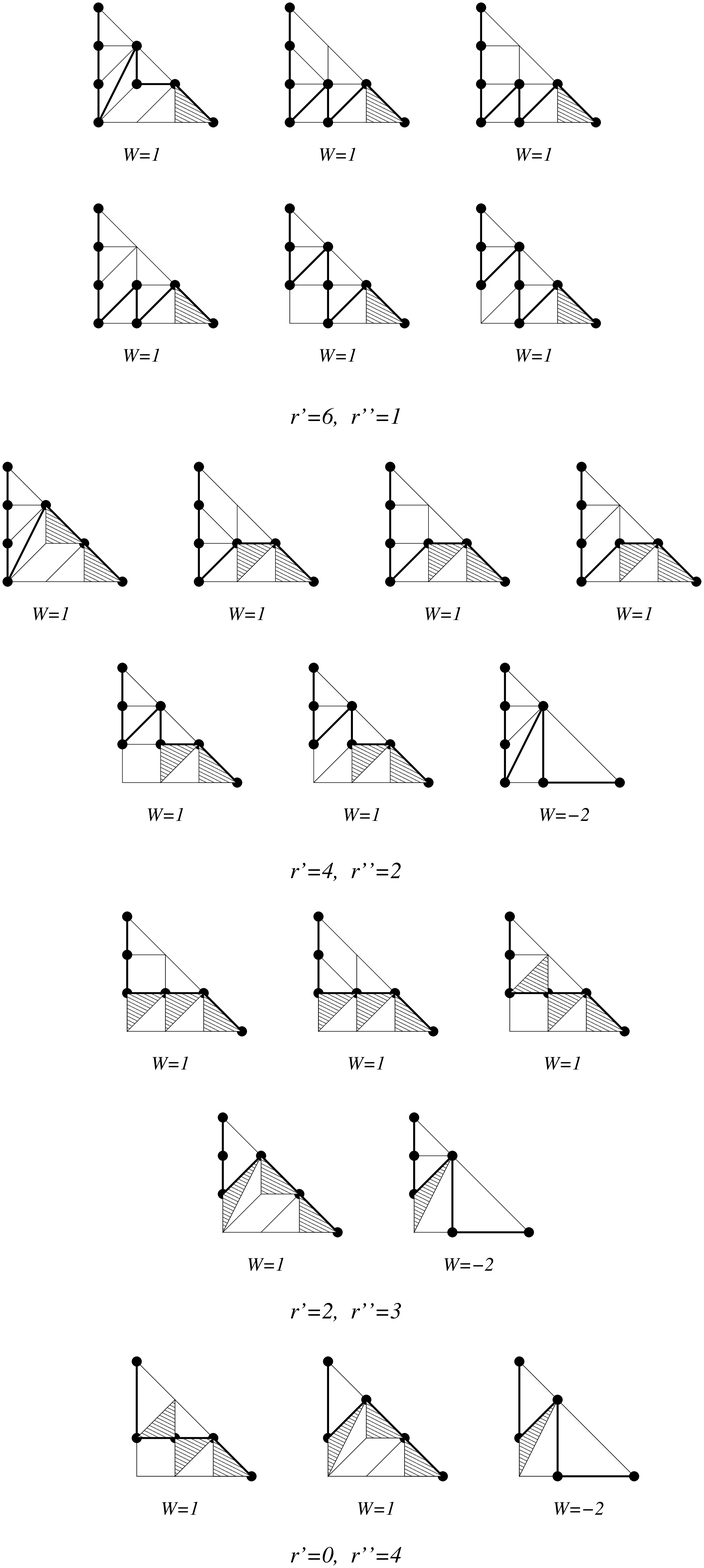}
\end{center}
\caption{Tropical count of $W_{r''}(\PP^2,3L)$} \label{fign6}
\end{figure}
Let the non-negative integers $r'$ and $r''$ satisfy $r'+2r''=8$.
Consider the pencil of real plane cubics passing through $r'$
generic real points and through $r''$ generic pairs of imaginary
conjugate points. Clearly, it has $r'+1$ real base points. Denote
by $n_+$ and $n_-$ the numbers of real cubics in the pencil, which
have a non-solitary real node or a solitary real node,
respectively. Blowing up the $r'+1$ real base points of the pencil
and integrating with respect to the Euler characteristic over the
blown-up real projective plane, we obtain (cf. \cite{Wel,Wel1})
\begin{equation}W_{r''}(\PP^2,3L)=n_+-n_-=r'=8-2r''\ .\label{enn50}\end{equation} Take a linear function
$\lam(x,y)=x-\bet y$ with $0<\bet\ll 1$, and put ${\mathcal
W}=\{1,...,r'\}$. Figure \ref{fign6} presents all consistent
subdivisions of the Newton triangle
$\Del=\conv\{((0,0),(0,3),(3,0)\}$, constructed along the
procedure of section \ref{secn10}, and their Welschinger numbers,
confirming formula (\ref{enn50}).

Notice that $W_4(\PP^2,3)=0$. In fact, there is a configuration of
four pairs of imaginary conjugate points in $\PP^2$, which defines
a pencil of non-singular real cubics. We present the following
example:
$F:=\alp(x_2^3+x_1^2x_2+x_0^2x_2)+\bet(x_0x_2^2+x_1^3+x_0^3)=0$.
By a routine computation, the system of equations for singular
points
$$\begin{cases}&F_{x_0}=3\alp x_0^2+\alp x_2^2+2\bet x_0x_2=0,\\
&F_{x_1}=3\alp x_1^2+2\bet x_1x_2=0,\\
&F_{x_2}=2\alp x_0x_2+\bet x_0^2+\bet x_1^2+3\bet
x_2^2=0\end{cases}$$ can be reduced to a relation for $\alp$ and
$\bet$, which has no real solutions except for $\bet=0$, which
leaves the only real singular curve $x_2^3+x_1^2x_2+x_0^2x_2=0$.
The latter curve has two imaginary nodes. Hence, by a small
variation of the pencil, we get rid of the real singular curves.

\bibliographystyle{amsplain}

\end{document}